\documentclass[11pt,a4paper,reqno]{amsart}
\usepackage[T1]{fontenc}
\usepackage{amsmath}
\usepackage{amsthm}
\usepackage{amsfonts}
\usepackage{amssymb}
\usepackage{graphicx}
\usepackage{amsbsy}
\newtheorem{theorem}{Theorem}[section]
\newtheorem{lemma}[theorem]{Lemma}
\newtheorem{ass}[theorem]{Assumption}
\newtheorem{proposition}[theorem]{Proposition}
\newtheorem{corollary}[theorem]{Corollary}

\theoremstyle{remark}
\newtheorem{remark}[theorem]{\it \bf{Remark}\/}
\newenvironment{acknowledgement}{\noindent{\bf Acknowledgement.~}}{}
\numberwithin{equation}{section}
\catcode`@=11
\def\section{\@startsection{section}{1}%
  \z@{1.5\linespacing\@plus\linespacing}{.5\linespacing}%
  {\normalfont\bfseries\large\centering}}
\catcode`@=12
\def\pa{\partial}
\def\na{\nabla}
\def\NN{\mathbb{N}}
\def\RR{\mathbb{R}}
\def\ds{\displaystyle}
\def\ni{\noindent}
\def\bs{\bigskip}
\def\ms{\medskip}

\def\eps{\varepsilon}
\def\fref#1{{\rm (\ref{#1})}}
\def\pref#1{{\rm \ref{#1}}}
\def\ints{\int_{\RR}}
\def\intd{\int_{\RR^{2}}}

\def\im{{\mathcal I} \hspace{-2pt}{\textit m}\,}

\def\calC{{\mathcal C}}

\def\calO{{\mathcal O}}

\def\calE{{\mathcal E}}

\def\Pe{\Psi^{\varepsilon}}

\def\Phe{\Phi^{\varepsilon}}

\def\Ve{V^{\varepsilon}}

\def\re{r^{\varepsilon}}

\catcode`@=11
\def\supess{\mathop{\operator@font Sup\,ess}}
\catcode`@=12
\newcommand{\be}{\begin{equation}}
\newcommand{\ee}{\end{equation}}
\newcommand{\bea}{\begin{eqnarray}}
\newcommand{\eea}{\end{eqnarray}}
\newcommand{\bee}{\begin{eqnarray*}}
\newcommand{\eee}{\end{eqnarray*}}
\title[2DEG in a strong magnetic field]{An effective mass theorem for the bidimensional electron gas in a strong magnetic field}
\author[F. Delebecque-Fendt]{Fanny Delebecque-Fendt}
\address{IRMAR, Universit\'e Rennes 1, France}
\email{fanny.fendt@univ-rennes1.fr}
\author[F. M\'ehats]{Florian M\'ehats}
\address{IRMAR, Universit\'e de Rennes 1, France}
\email{florian.mehats@univ-rennes1.fr}

\begin{document}

\begin{abstract}
We study the limiting behavior of a singularly perturbed Schr\"odinger-Poisson system describing a 3-dimensional electron gas  strongly confined in the vicinity of a plane $(x,y)$ and subject to a strong uniform magnetic field in the plane of the gas. The coupled effects of the confinement and of the magnetic field induce fast oscillations in time that need to be averaged out. We obtain at the limit a system of 2-dimensional Schr\"odinger equations in the plane $(x,y)$, coupled through an effective selfconsistent electrical potential. In the direction perpendicular to the magnetic field, the electron mass is modified by the field, as the result of an averaging of the cyclotron motion. The main tools of the analysis are the adaptation of the second order long-time averaging theory of ODEs to our PDEs context, and the use of a Sobolev scale adapted to the confinement operator.
\end{abstract}


\maketitle


\sloppy
\section{Introduction}
\subsection{The singularly perturbed problem}
Many electronic devices are based on the quantum transport of a bidimensional electron gas (2DEG) artificially confined in heterostructures at nanometer scales, see e.g. \cite{AFS,Bas,FG,VW}. In this article, we derive an asymptotic model for the quantum transport of a 2DEG subject to a strong uniform magnetic field which is parallel to the plane of the gas. The aim of this paper is to understand how the cyclotron motion competes with the effects of the potential confining the electrons and the nonlinear effects of the selfconsistent Poisson potential. Our tool is an asymptotic analysis from a singularly perturbed Schr\"odinger-Poisson system towards a reduced model of bidimensional quantum transport. In particular, we generalize in this context the notion of {\em cyclotron effective mass}, usually explicitely calculated in the simplified situation of a harmonic confinement potential \cite{FG,SJ}.

Our starting model is thus the 3D Schr\"odinger-Poisson system, singularly perturbed by a confinement potential and the strong magnetic field. The three-dimensional space variables are denoted by $(x,y,z)$ and the associated canonical basis of $\RR^3$ is denoted by $(e_x,e_y,e_z)$. The particles are subject to three effects: a confinement potential depending on the $z$ variable, a uniform magnetic field applied to the gas along the $e_y$ axis, and the selfconsistent Poisson potential. Given a small parameter $\eps>0$, which is the typical extension of the 2DEG in the $z$ direction, our starting model is the following dimensionless Schr\"odinger-Poisson system:
\begin{equation}
\label{schrod}
i\partial_{t}\Pe=\frac{1}{\eps^2}\left(-\partial_{z}^2+B^2z^2+V_c(z)\right)\Pe-\frac{1}{\eps}2iBz\pa_x\Pe-\Delta_{x,y}\Pe+\Ve \Pe\,,
\end{equation}
\begin{equation}
\label{cauchy}
\Pe(0,x,y,z)=\Psi_0(x,y,z),
\end{equation}
\begin{equation}
\label{poisson}
\Ve(t,x,y,z)=\frac{1}{4\pi \re}\ast |\Pe|^2,
\end{equation}   
where we have denoted
\be
\label{re}
\re(x,y,z)=\sqrt{x^2+y^2+\varepsilon^2z^2}.
\ee
The scaling is discussed in the next subsection. This system describes the transport of electrons under the action of:
\begin{itemize}
\item[--] The applied confinement potential $\frac{1}{\eps^2}V_c(z)$, nonnegative, such that $V_c(z)\to +\infty$ as $|z|\to +\infty$. The precise assumptions of this potential are made below in Assumptions \ref{confinement} and \ref{ass2}. 
\item[--] The applied uniform magnetic field $\frac{B}{\eps}e_y$ (with $B>0$ fixed), which derives from the magnetic potential $\frac{1}{\eps}Bze_x$. We have chosen to work in the Landau gauge.
\item[--] The Poisson selfconsistent potential $V^\eps$.
\end{itemize}
Note that \fref{schrod} is equivalent to
\begin{equation}
\label{schrodA}
i\partial_{t}\Pe=\frac{1}{\eps^2}\left(-\partial_{z}^2+V_c(z)\right)\Pe+\left(i\pa_x-\frac{Bz}{\eps}\right)^2\Pe-\pa^2_y\Pe+\Ve\Pe\,.
\end{equation}
The goal of this work is to exhibit an asymptotic system for \fref{schrod}, \fref{cauchy}, \fref{poisson} as $\eps\to 0$.

Let us end this subsection with short bibliographical notes. In a linear setting, quantum motion constraint on a manifold has been studied for a long time by several authors, see \cite{dC,DE,FH,T} and references therein. Nonlinear situations were studied more recently. The approximation of the Schr\"odinger-Poisson system with no magnetic field was studied when the electron gas is constraint in the vicinity of a plane in \cite{BAMP,P} and when the gas is constraint on a line in \cite{BCFM}. When the nonlinearity depends locally on the density, as for the Gross-Pitaevskii equation, asymptotic models for confined quantum systems were studied in \cite{BMSW,BACM,CMS}. In classical setting, collisional models in situations of strong confinement have been studied in \cite{DPV}. Finally, let us draw a parallel with the problem of homogenization of the Schr\"odinger equation in a large periodic potential, studied in \cite{allaire} and \cite{sparber}. At the limit $\eps\to 0$, as noted above, we will obtain an homogenized system which takes the form of bidimensional Schr\"odinger equations with an effective mass in the $x$ direction. However, this phenomenon is due to an averaging of the cyclotron motion induced by a strong magnetic field, and is not exactly the same notion as the usual effective mass for the transport in a lattice or in a crystal. Nevertheless, it is interesting to observe that the scaling used in \cite{allaire, sparber} in the case of a strong periodic potential is similar to the strong confinement scaling used in the present paper.

\subsection{The physical scaling}

In order to clarify the physical assumptions underlying our singularly perturbed system, let us derive \fref{schrod}, \fref{cauchy}, \fref{poisson} from the Schr\"odinger-Poisson system written in physical variables. This system reads as follows:
\begin{equation}
 \label{schrodphys}
 i \hbar\partial_{\mathbf t}{\mathbf \Psi}=\frac{1}{2m}\left(i\hbar\nabla -\frac{e{\mathbf B}}{c}{\mathbf z }e_{{\mathbf x}}\right)^2 {\mathbf \Psi}+e{\mathbf {V_{c}}}{\mathbf \Psi}+e{\mathbf V} {\mathbf \Psi},
 \end{equation}
\begin{equation}
\label{poissonphys}
{\mathbf V}=\frac{e}{4\pi \epsilon\sqrt{\mathbf x^2+\mathbf y^2+\mathbf z^2}}\ast\left(|{\mathbf \Psi}|^2\right).
\end{equation}
Each dimensionless quantity in \fref{schrod}, \fref{cauchy}, \fref{poisson} is the associated physical quantity normalized by a typical scale:
\be
\label{scale1}x=\frac{\mathbf x}{\overline{x}},\;y=\frac{\mathbf y}{\overline{y}},\;z=\frac{\mathbf z}{\overline{z}},\;|\Psi^\eps|^2=\frac{ |{\mathbf \Psi}|^2}{\overline N},\;V_{c}=\frac{\mathbf{V_{c}}}{\overline{V_{c}}},\;V^\eps=\frac{\mathbf{V}}{\overline{V}},\;B=\frac{\mathbf B}{\overline B}.
\ee
Now we introduce two energy scales in this problem: a strong energy $E_{conf}$, which will be the energy of the confinement in $z$ and of the magnetic effects, and a transport energy $E_{transp}$, which will be the typical energy of the longitudinal transport in $(x,y)$ and also of the selfconsistent effects. We introduce the following small dimensionless parameter:
\be
\label{assenergy}
\eps=\left(\frac{E_{transp}}{E_{conf}}\right)^{1/2}\ll 1.
\ee
Then our scaling assumptions are the following. We set to the scale $E_{conf}$ the confinement potential, the magnetic energy and the kinetic energy along $z$:
\be\label{scale2}
E_{{conf}}:=e\overline{V_{c}}=\frac{1}{2}m\left(\frac{e\overline B}{mc}\right)^2\overline z^2=\frac{\hbar^2}{2m\overline z^2}
\ee
and we set to the scale $E_{transp}$ the selfconsistent potential energy, the kinetic energies along $x$ and $y$ and we finally choose a time scale adapted to this energy:
\be
\label{scale3}E_{transp}:= e{\overline V}= \frac{e^2\overline N \,\overline x \,\overline z}{\epsilon}=\frac{\hbar^2}{2m \overline x^2}=\frac{\hbar^2}{2m \overline y^2}= \frac{\hbar}{\overline t}.
\ee
By inserting \fref{scale1} in \fref{schrodphys}, \fref{poissonphys}, then by using \fref{assenergy}, \fref{scale2} and \fref{scale3}, we obtain directly our singularly pertubed problem \fref{schrod}, \fref{poisson}. Note that \fref{scale2} and \fref{scale3} imply that $\eps$ is also the ratio between the transversal and the longitudinal space scales:
$$\eps=\frac{\overline z}{\overline x}=\frac{\overline z}{\overline y}.$$

\subsection{Heuristics in a simplified case}
\label{heuristics}
In this section, we analyze a very simplified situation where analytic calculations can be directly done. We assume here that $V_c$ is a harmonic confinement potential and we neglect the Poisson potential $V^\eps$. We formally analyze the heuristics in this simplified case, that will be further compared to our result obtained in the general case.

We thus consider here a new system, similar to \fref{schrod} where we prescribe $V_c(z)=\alpha^2z^2$, $\alpha>0$ and where the Poisson potential $V^\eps$ is replaced by 0:
\begin{equation}
\label{schrodpart}
i\partial_{t}\Pe=\frac{1}{\eps^2}\left(-\partial_{z}^2+(\alpha^2+B^2)z^2\right)\Pe-\frac{1}{\eps}2iBz\pa_x\Pe-\Delta_{x,y}\Pe\,,
\end{equation}
\begin{equation}
\label{cauchypart}
\Pe(0,x,y,z)=\Psi_0(x,y,z).
\end{equation}
In this situation, there is a trick which enables to transform the equation. Indeed, by remarking that
\begin{eqnarray*}
&&-\partial_{z}^2+(\alpha^2+B^2)z^2-2iB\eps z\partial_x-\eps^2\partial^2_x\\
&&\qquad=-\partial_{z}^2+(\alpha^2+B^2)\left(z-\frac{B}{\alpha^2+B^2}\,i\eps\partial_x\right)^2-\frac{\alpha^2}{\alpha^2+B^2}\eps^2\partial^2_x\,,
\end{eqnarray*}
we obtain that \fref{schrodpart} is equivalent to
\be
i\partial_{t}\Pe=\frac{1}{\eps^2}\left[-\partial_{z}^2+(\alpha^2+B^2)\left(z-\frac{B}{\alpha^2+B^2}\,i\eps\partial_x\right)^2\right]\Pe -\frac{\alpha^2}{\alpha^2+B^2}\partial^2_x\Pe-\pa_y^2\Pe\,.\label{schrodpart2}
\ee
Introduce now the following operator: for a function $u\in L^2(\RR^3)$, we set
$$(\Theta^\eps u)(x,y,z)=\mathcal F_x^{-1}\left(\mathcal F_x u (\xi,y,z+\frac{B}{\alpha^2+B^2}\,\eps\xi)\right),$$
where $\mathcal F_x$ denotes the Fourier transform in the $x$ variable. Note that this operator $\Theta^\eps$ is unitary on $L^2(\RR^3)$ and commutes with $\pa_x$ and $\pa_y$. Hence, we deduce from \fref{schrodpart2} and by direct calculations that the function $u^\eps=\Theta^\eps \Psi^\eps$ satisfies the following system:
\begin{equation*}
i \partial_{t}u^{\eps}=\frac{1}{\varepsilon^2}\widetilde H_{z}u^{\eps}-\frac{\alpha^2}{\alpha^2+B^2}\partial^2_xu^{\eps}-\partial^2_yu^{\eps},\quad u^\eps(t=0)=\Theta^\eps \Psi_0\,,
\end{equation*}
where
$$\widetilde H_{z}=-\partial_{z}^2+(\alpha^2+B^2)z^2.$$
Let us now filter out the oscillations by introducing the new unknown
$$\Phi^\eps=\exp(it\widetilde H_z/\eps^2)u^\eps.$$
Again, the operator $\exp(it\widetilde H_z/\eps^2)$ commutes with $\pa_x$, $\pa_y$ and, finally, the following equation is equivalent to \fref{schrodpart}:
\begin{equation}
\label{phipart}
i \partial_{t}\Phi^{\eps}=-\frac{\alpha^2}{\alpha^2+B^2}\partial^2_x\Phi^{\eps}-\partial^2_y\Phi^{\eps},\qquad \Phi^\eps(t=0)=\Theta^\eps \Psi_0\,.
\end{equation}
As $\eps \to 0$, it is not difficult to see that, for sufficiently smooth initial data, we have $\Theta^\eps \Psi_0\to \Psi_0$. Therefore, one can show that, in adapted functional spaces, we have $\Phi^\eps \to \Phi$ as $\eps\to 0$, with $\Phi$ solution of the limit system:
\begin{equation}
\label{schrodlimitpart}
i \partial_{t}\Phi=-\frac{\alpha^2}{\alpha^2+B^2}\partial^2_x\Phi-\partial^2_y\Phi,\qquad \Phi(t=0)= \Psi_0\,.
\end{equation}
This equation is a bidimensional Schr\"odinger equation with an anisotropic operator that can be interpreted as follows. Whereas, as expected, the dynamics in the $y$ is not perturbed by the magnetic field (since it is parallel to $y$), in the $x$ direction the electrons are transported as if their mass was augmented by a factor $\frac{\alpha^2+B^2}{\alpha^2}>1$. This coefficient is called the (dimensionless) {\em electron cyclotron mass} \cite{FG,SJ}.

In this article, the model that we want to treat is the nonlinear system \fref{schrod}, \fref{poisson}, with a general confinement potential $V_c$ instead of $\alpha^2 z^2$ and the selfconsistent Poisson potential. Consequently, it is not possible to simplify the equation \fref{schrod} by the above trick. Moreover, the potential $V^\eps$ depends on the $z$ variable and on the function $\Pe$ itself. Therefore, one has to be careful for instance when filtering out the fast oscillations by applying the operator $\exp(it\widetilde H_z/\eps^2)$, since in this nonlinear framework some interference effects between the elementary waves might appear. In this article, we present a general strategy that enables to overcome these difficulties. The strategy will be inspired from \cite{BACM} where the nonlinear Schr\"odinger equation under strong partial confinement was analyzed. Two main differences appear here. First, the Poisson nonlinearity is nonlocal, which requires specific estimates. Observe that, at the limit $\eps\rightarrow 0$, the nonlinearity in the present paper reads $\frac{1}{4\pi |x|}\ast\int |\psi|^2dz$ and does not depend on $z$. This makes an important difference with the case of \cite{BACM}, in particular no resonance effects due to the nonlinearity will appear. Second, the magnetic field induces in \fref{schrod} a singular term at an intermediate scale $\frac{1}{\eps}$ between the confinement operator (at the scale $\frac{1}{\eps^2}$) and the nonlinearity (at the scale $\frac{1}{\eps^0}$). Hence, compared to \cite{BACM}, the average techniques have to be pushed to the order two and resonance effects will finally appear here due to this magnetic term.

\subsection{Main result}

Consider the system \fref{schrod}, \fref{cauchy}, \fref{poisson}. We assume that the confinement potential $V_c$ satisfies two assumptions. The first one concerns the behavior of this function at the infinity.
\begin{ass}
\label{confinement}
The potential $V_c$ is a $\calC^\infty$ nonnegative {\em even} function such that
\be
\label{H1}
a^2|z|^2\leq V_c(z)\leq C|z|^M \quad\mbox{for }|z|\geq 1,
\ee
where $a>0$, $M>0$, and
\be
\label{H2}
\frac{|\partial_{z}V_{c}(z)|}{V_c(z)}=\calO\left(|z|^{-M'}\right),
\quad \frac{|\partial^k_zV_c(z)|}{V_c(z)}=\calO(1)\mbox{ for all }k\in \NN^*,
\ee
as $|z|\to +\infty$, where $M'>0$.
\end{ass}
Note that a smooth even potential of the form $V_c(z)=C |z|^s$ for $|z|\geq |z_0|$, with $C>0$, $s\geq 2$, satisfies these assumptions. In particular the harmonic potential $V_{c}=a^2 z^2$ fits these conditions.

Let us discuss on the assumptions. The assumption that the function $V_c(z)$ is even is important in our analysis, see e.g. Step 4 in subsection \ref{sketch}. The left inequality in the first condition \fref{H1} implies that $V_c$ tends to $+\infty$ as $|z|\to +\infty$. The fact that $V_c(z)\geq a^2 z^2$ is not essential in our analysis but simplifies it (see below, it allows to give a simple characterization of the energy space related to our system). As it is well-known \cite{RS}, the spectrum of operator $H_z$ defined by
\be
\label{Hz}
H_z=-\pa^2_z+B^2z^2+V_c(z).
\ee
is discrete, when $H_z$ is considered as a linear, unbounded operator over $L^2(\RR)$, with domain 
$$D(H_z)=\{u\in L^2(\RR),\,\, H_zu\in L^2(\RR)\}.$$
The complete sequence of eigenvalues of $H_z$ will be denoted by $(E_{p})_{p\in \NN}$, taken strictly increasing with $p$ (recall indeed that in dimension 1 the eigenvalues are simple), and the associated Hilbert basis of real-valued eigenfunctions will be denoted by $(\chi_{p}(z))_{p\in \NN}$. The right inequality in \fref{H1} and the second condition \fref{H2} are more technical and are here to simplify the use of a Sobolev scale based on the operator $H_z$, which is well adapted to our problem. More precisely, these assumptions are used in Lemma \ref{lemBACM}.

The second assumption on $V_c$ concerns the spectrum of the confinement operator $H_z$.
\begin{ass}
\label{ass2}
The eigenvalues of the operator $H_z$ defined by \fref{Hz} satisfy the following property: there exists $C>0$ and $n_0\in \NN$ such that
$$\forall p\in \NN,\quad E_{p+1}-E_p\geq C(1+p)^{-n_0}.$$
\end{ass}
The most simple situation where \fref{ass2} is satisfied is when there exists a uniform gap between the eigenvalues: for all $p\in \NN^*$, $ E_{p+1}-E_p\geq C_0>0$. Note that in this case we have $n_0=0$. This property is true in the following examples.
\begin{itemize}
\item[--] If $V_c(z)=a^2z^2+V_1(z)$, with $\|V_1\|_{L^\infty}<2\sqrt{a^2+B^2}$. Indeed, in this case the perturbation theory gives $|E_p-(2p+1)\sqrt{a^2+B^2}|<\|V_1\|_{L^\infty}$.
\item[--] If $V_c(z)\sim a|z|^s$ as $|z|\to +\infty$, with $s>2$. Indeed, in this case the Weyl asymptotics \cite{egorov} gives $E_p\sim Cp^{\frac{2s}{s+2}}$,  so $E_{p+1}-E_p\to +\infty$ as $p\to +\infty$.
\end{itemize}

Let us now give a few indications on the Cauchy problem for \fref{schrod}, \fref{cauchy}, \fref{poisson}. This system benefits from two conservation laws, the mass and energy conservations:
\be
\label{cons}
\forall t\geq 0,\quad \|\Psi^\eps(t)\|_{L^2}^2=\|\Psi_0\|_{L^2}^2,\qquad \calE(\Pe(t))=\calE(\Psi_0),
\ee
where the total energy of the wavefunction $\Psi^\eps$ is defined by
\bea
&\calE(\Psi^\eps)=&\frac{1}{\eps^2}\|\pa_z \Psi^\eps\|^2_{L^2}+\frac{1}{\eps^2}\|\sqrt{V_c}\Psi^\eps\|_{L^2}^2+\frac{1}{\eps^2}\|(\eps\pa_x+iBz)\Psi^\eps\|_{L^2}^2\nonumber\\
&&+\| \pa_y \Psi^\eps\|_{L^2}^2+\frac{1}{2}\|\sqrt{V^\eps}\Psi^\eps\|_{L^2}^2\,.\label{energy}
\eea
For fixed $\eps>0$, the Cauchy theory for the Schr\"odinger-Poisson with a constant uniform magnetic field was solved in \cite{cazenave-esteban,debouard} in the energy space. It is not difficult to adapt these proofs (see also the reference book \cite{CAZ}) to our case where an additional confinement potential is applied. The energy space in our situation is the set of functions  $u$ such that $\calE(u)$ is finite:
\bee
B^1=\left\{u\in L^2(\RR^3):\,\pa_zu\in L^2(\RR^3),\, \sqrt{V_c}u\in L^2(\RR^3),\, \pa_yu\in L^2(\RR^3)\right.\qquad &&\\
\mbox{and } \left.\left(\pa_x+\frac{iBz}{\eps}\right)u\in L^2(\RR^3)\right\}.&&
\eee
This space seems to depend on $\eps$, which would not be convenient for our asymptotic analysis. In fact, it does not. Indeed, thanks to our assumption \fref{H1} on the confinement potential, one has
$$\|zu\|_{L^2}\leq \frac{1}{a}\|\sqrt{V_c}u\|_{L^2},$$
so $u\in B^1$ implies that $zu\in L^2$ and thus $\pa_x u\in L^2$. Hence one has
$$B^1=\left\{u\in H^1(\RR^3): \sqrt{V_c}u\in L^2(\RR^3)\right\}$$
and, on this space, we will use the following norm independent of $\eps$:
\bea
&\|u\|_{B^1}^2&=\|(I-\Delta_{x,y}+H_z)^{1/2}u\|_{L^2}^2\nonumber\\
&&=\|u\|_{L^2}^2+\|(-\Delta_{x,y})^{1/2} u\|_{L^2}^2+\|(H_z)^{1/2} u\|_{L^2}^2\nonumber\\
&& =\|u\|_{H^1}^2+\|\sqrt{V_c}u\|_{L^2}^2+B^2\|zu\|_{L^2}^2,\label{equivB11}
\eea
where we used the selfadjointness and the positivity of $-\Delta_{x,y}$ and of the operator $H_z$ defined by \fref{Hz}, and where $I$ denotes the identity operator. In this paper, we will assume that the initial datum $\Psi_0$ in \fref{cauchy} belongs to this space $B^1$. Then, for all $\eps>0$,  the system \fref{schrod}, \fref{cauchy}, \fref{poisson} admits a unique global solution $\Pe\in C^0([0,+\infty), B^1)$. Our aim is to analyze the asymptotic behavior of $\Pe$ as $\eps\to 0$.

\bs
We are now in position to state our main results. Here and throughout this paper, we will use the notation
\be
\label{langle}
\forall u\in L^1_z(\RR),\qquad \langle u\rangle =\int_{\RR}u(z)\,dz.
\ee
Let us introduce the limit system. First define the following coefficients
\be
\label{alphap}
\forall p\in \NN,\qquad \alpha_p=1-\sum_{q\neq p}\frac{\left\langle 2Bz \chi_p\chi_q\right\rangle^2}{E_q-E_p}\,,
\ee
where we recall that $(E_p,\chi_p)_{p\in\NN}$ is the complete sequence of eigenvalues and eigenfunctions of the operator $H_z$ defined by \fref{Hz}. Then, we introduce the following infinite dimensional, nonlinear and coupled differential system on the functions $\phi_p(t,x,y)$:
\be
\label{limitS}
\forall p\in \NN,\qquad i\pa_t \phi_p=-\alpha_p\,\pa^2_x \phi_p-\pa^2_y\phi_p+W\phi_p\,,\quad \phi_p(t=0)=\left\langle\Psi_0\,\chi_p\right\rangle\,,
\ee
\begin{equation}
\label{limitP}
W=\frac{1}{4\pi \sqrt{x^2+y^2}}\ast\left(\sum_{p\in \NN}|\phi_{p}|^2\right).
\end{equation}
Note that the convolution in \fref{limitP} holds on the variables $(x,y)\in \RR^2$. The equation \fref{limitP} is nothing but the Poisson equation for a measure valued distribution of mass whose support is constrained to the plane $z=0$:
$$W(t,x,y)=\left.\left[\frac{1}{4\pi \sqrt{x^2+y^2+z^2}}\ast\left(\sum_{p\in \NN}|\phi_{p}(t,x,y)|^2\delta_{z=0}\right)\right]\right|_{z=0}.$$
In order to compare with $\Pe$, we introduce the following functions:
\be
\label{psilimit}
\Phi(t,x,y,z)=\sum_{p\in \NN}\phi_p(t,x,y)\,\chi_p(z), \quad \Psi^{\eps}_{app}(t,x,y,z)=\sum_{p\in \NN}e^{-itE_p /\eps^2}\,\phi_p(t,x,y)\,\chi_p(z).
\ee
Remark that $\Psi^{\eps}_{app}$ can be deduced from $\Phi$ through the application of the operator $e^{itH_z/\eps^2}$, unitary on $B^1$:
$$\Psi^{\eps}_{app}=e^{-itH_z/\eps^2}\,\Phi.$$ This explicit relation is the only dependency in $\eps$ of the limit system \fref{limitS}, \fref{limitP}, \fref{psilimit}. Our main result is the following theorem.
\begin{theorem}
\label{theo}
Assume that $V_c$ satisfies Assumptions \pref{confinement} and \pref{ass2} and let $\Psi_0\in B^1$. For all $\eps\in(0,1]$, denote by $\Pe\in C^0([0,+\infty), B^1)$ the unique global solution of the initial system \fref{schrod}, \fref{cauchy}, \fref{poisson}. Then the following holds true.\\
(i) The limit system \fref{limitS}, \fref{limitP}, \fref{psilimit} admits a unique maximal solution $\Psi^{\eps}_{app}\in C^0([0,T_{max}),B^1)$, where $T_{max}\in (0,+\infty]$ is independent of $\eps$. If $T_{max}<+\infty$ then $\|\Psi^{\eps}_{app}(t,\cdot)\|_{B^1}\to +\infty$ as $t\to T_{max}$.\\
(ii) For all $T\in (0,T_{max})$, we have
  $$ \lim_{\eps\to 0}\left\|\Psi^\eps -\Psi^{\eps}_{app}\right\|_{C^0([0,T],B^1)}=0.
$$
\end{theorem}

\bs
\ni
{\bf Comments on Theorem \pref{theo}.}

\ms
\ni
{\em 1. The cyclotron effective mass.} Theorem \ref{theo} thus states that, on all time intervals where the limit system \fref{psilimit}, \fref{limitS}, \fref{limitP} is well-posed, the solution $\Psi^\eps$ of the singularly perturbed system \fref{schrod}, \fref{cauchy}, \fref{poisson} is close to $\Psi^{\eps}_{app}$. As expected, the dynamics in the $y$ direction, ie parallel to the magnetic field, is not affected by the magnetic field, since the operator is still $-\pa_y^2$. On the other hand, the situation is different in the direction $x$ and the averaging of the cyclotron motion results in a multiplication of the operator $-\pa_x^2$ by the factor $\alpha_p$ which only depends on $V_c$ and $B$. The coefficient $\frac{1}{\alpha_p}$ plays in \fref{limitS} the role of an effective mass in the direction perpendicular to the magnetic field. We find that the effective mass in the Schr\"odinger equation for the mode $p$ depends on the index $p$ of this mode. We do not know whether these coefficients are positive for a general $V_c$. 

Notice that the effective mass could be predicted heuristically by the following argument. Denoting by $k_x$, $k_y$ the wavevectors of the 2DEG in the plane $(x,y)$, the electron dispersion relation $E_p(k_x,k_y)$ in the transversal subbands can be written from \fref{schrod} by computing the eigenvalues of the operator
$$\frac{1}{\eps^2}\left(-\frac{d^2}{dz^2}+B^2z^2+V_c(z)+2\eps Bzk_x+\eps^2 k_x^2+\eps^2k_y^2\right).$$
Since $\eps$ is small, an approximation of $E_p(k_x,k_y)$ can be computed thanks to perturbation theory, which gives the following parabolic band approximation:
$$E_p(k_x,k_y)=\frac{E_p}{\eps^2}-k_x^2 \sum_{q\neq p}\frac{\left\langle 2Bz \chi_p\chi_q\right\rangle^2}{E_q-E_p}+k_x^2+k_y^2+o(1).$$
We can read on this formula that the effective mass is 1 in the $y$ direction and is $\alpha_p^{-1}$ according to \fref{alphap} in the $x$ direction. Note that the specific case of the harmonic potential is treated below (see comment $3$).

\ms
\ni
{\em 2. Conservation of the energy for the limit system.} Let us write the conservation of the energy for the limit system. The total energy for this system can be splitted into a confinement energy $\calE_{conf}(\Phi)$ and a transport energy $\calE_{tr}(\Phi)$ defined by
\bea
&&\hskip -1cm \calE_{conf}(\Phi)=\sum_{p\in \NN}E_p\,\|\phi_p\|_{L^2}^2\,,\label{econf}\\
\nonumber&&\hskip -1cm \calE_{tr}(\Phi)=\sum_{p\in \NN}\alpha_p \|\pa_x \phi_p\|_{L^2}^2+\sum_{p\in \NN}\|\pa_y \phi_p\|_{L^2}^2\\
&&\hskip -0.5cm +\frac{1}{2}\sum_{p,q}\int_{\RR^4}\frac{1}{4\pi \sqrt{|x-x'|^2+|y-y'|^2}}|\phi_p(x,y)|^2|\phi_q(x',y')|^2\,dxdydx'dy'.\label{etransp}
\eea
An interesting property is that these two quantities are separately conserved by the limit system. If $\Psi^{\eps}_{app}$ solves  \fref{limitS}, \fref{limitP}, \fref{psilimit}, then, for all $t\in [0,T]$, we have
\be
\label{consenerlimit}
\calE_{conf}(\Psi^{\eps}_{app}(t))=\calE_{conf}(\Psi^{\eps}_{app}(0))\quad \mbox{and}\quad \calE_{tr}(\Psi^{\eps}_{app}(t))=\calE_{tr}(\Psi^{\eps}_{app}(0)).
\ee
In particular, by summing up the two equalities in \fref{consenerlimit}, we obtain the following conservation property:
\be
\label{enermod}
\calE_{conf}(\Psi^{\eps}_{app}(t))+\calE_{tr}(\Psi^{\eps}_{app}(t))=\calE_{conf}(\Psi^{\eps}_{app}(0))+\calE_{tr}(\Psi^{\eps}_{app}(0)).
\ee
Note that, in the general case, we do not know whether the energy defined by \fref{etransp} is the sum of nonnegative terms. This point is related to the fact that the well-posedness for $t\in[0,+\infty)$ of the Cauchy problem for the nonlinear system \fref{limitS}, \fref{limitP} is an open issue. Nevertheless, when the $\alpha_p$ are such that the energy is coercive on $B^1$, ie when we have
 \be
 \label{coercive}
 \forall \Phi\in B^1,\quad C_0\|\Phi\|_{B^1}^2\leq \calE_{conf}(\Phi)+\calE_{tr}(\Phi)\leq C_1\|\Phi\|_{B^1}^2+C_2\|\Phi\|_{B^1}^4,
 \ee
 with a constant $C_0>0$ independent of $\eps$, then the maximal solution of \fref{limitS}, \fref{limitP} is globally defined: $T_{max}=+\infty$.
 \begin{corollary}[Global in time convergence]
\label{theo2}
Under the assumptions of Theorem \pref{theo},  assume moreover that there exists $0<\underline \alpha<\overline \alpha$ such that the coefficients $\alpha_p$ defined by \fref{alphap} satisfy the following condition:
\be
\label{delta}
\forall p\in \NN,\qquad \underline\alpha\leq\alpha_p\leq \overline \alpha.
\ee
Then the system \fref{psilimit}, \fref{limitS}, \fref{limitP} admits a unique global solution $\Psi^{\eps}_{app}\in C^0([0,+\infty),B^1)$ and, for all $T>0$, we have  $$ \lim_{\eps\to 0}\left\|\Psi^\eps -\Psi^{\eps}_{app}\right\|_{C^0([0,T],B^1)}=0,
$$
where $\Psi^\eps\in C^0([0,+\infty),B^1)$ denotes the solution of \fref{schrod}, \fref{cauchy}, \fref{poisson}. 
\end{corollary}
The proof of this corollary is immediate and will not be detailed in this paper. Indeed, remarking that \fref{delta} implies \fref{coercive}, we obtain that the solution $\Psi^{\eps}_{app}(t)$ of \fref{limitS}, \fref{limitP} satisfies the following uniform bound:
 $$
 \|\Psi^{\eps}_{app}(t)\|_{B^1}^2\leq  C\left(\widetilde\calE_{conf}(\Psi^{\eps}_{app}(t))+\widetilde\calE_{tr}(\Psi^{\eps}_{app}(t))\right)= C\left(\widetilde\calE_{conf}(\Psi_0)+\widetilde\calE_{tr}(\Phi_0)\right),
 $$
 where the quantity in the right-hand side is finite as soon as $\Psi_0\in B^1$.
 
 \ms
 \ni
{\em 3. Case of harmonic confinement.} In the special case of a harmonic confinement potential $V_c(z)=a^2 z^2$, the eigenvalues and eigenfunctions of $H_z=-\pa_z^2+(a^2+B^2)z^2$ can be computed explicitely and one has
$$E_p=(2p+1)\sqrt{a^2+B^2},\qquad \chi_p(z)=(a^2+B^2)^{1/8}\, u_p\left((a^2+B^2)^{1/4}z\right),$$
where $(u_p)_{p\in \NN}$ are the normalized Hermite functions defined e.g. in \cite{M}, B 8 and satisfying $-u''+z^2u_p=(2p+1)u_p$.
The properties of the Hermite functions give
$$2z\chi_p=\frac{\sqrt{2(p+1)}}{(a^2+B^2)^{1/4}}\chi_{p+1}+\frac{\sqrt{2p}}{(a^2+B^2)^{1/4}}\chi_{p-1}\,,$$
and one can compute explicitely the coefficients
$$\alpha_p=1-B^2\frac{\langle 2z\chi_p\chi_{p+1}\rangle^2}{E_{p+1}-E_p}+B^2\frac{\langle 2z\chi_p\chi_{p-1}\rangle^2}{E_{p}-E_{p-1}}=\frac{a^2}{a^2+B^2}.$$
We thus recover here the coefficient found in subsection \ref{heuristics} in the simplified situation. Note that, in this case, condition \fref{delta} is satisfied and the convergence result holds on an arbitrary time interval. It is reasonable to conjecture that this condition \fref{delta} holds again when $V_c(z)=a^2z^2+V_1(z)$, where $V_1$ is a small perturbation.

 \ms
 \ni
{\em 4. Towards a more realistic model.}
Since we aim at describing the transport of electrons, which are fermions, our model should not be restricted to a pure quantum state. The following model describes the transport of an electron gas in a mixed quantum state and is more realistic:
\begin{equation}
\label{schrodmixed}
i\partial_{t}\Pe_j=\frac{1}{\eps^2}\left(-\partial_{z}^2+B^2z^2+V_c(z)\right)\Pe_j-\frac{1}{\eps}2iBz\pa_x\Pe_j-\Delta_{x,y}\Pe_j+\Ve\Pe_j\,,\quad \forall j,
\end{equation}
\begin{equation}
\label{cauchymixed}
\Pe_j(0,x,y,z)=\Psi_{j,0}(x,y,z),\quad \forall j,
\end{equation}
\begin{equation}
\label{poissonmixed}
\Ve(t,x,z)=\frac{1}{4\pi \re}\ast \rho^\eps,\qquad \rho^\eps=\sum_j \lambda_j |\Pe_j|^2,
\end{equation}   
where $\lambda_j$, the occupation factor of the state $\Pe_j$, takes into account the statistics of the electron ensemble and is fixed once for all at the initial time. Note that the Schr\"odinger equations \fref{schrodmixed} are only coupled through the selfconsistent Poisson potential. Therefore, we claim that our main Theorem \ref{theo}, which has been given for the sake of simplicity in the case of pure quantum state, can be extended to this system \fref{schrodmixed}, \fref{cauchymixed}, \fref{poissonmixed}, with appropriate assumptions on the initial data $(\Psi_{j,0})$\,.

Similarly, a given smooth external potential could be incorporated in the initial system. We also claim that our result can be easily adapted if we add in the right-hand side of \fref{schrod} a term of the form $V_{ext}(t,x,y,\eps z)\Psi^\eps$ (which is coherent with our scaling), and the result does not change qualitatively.

\subsection{Scheme of the proof}
\label{sketch}

In this section, we sketch the main steps of the proof of the main theorem.

\bs
\ni
{\em Step 1: a priori estimates}.

\ms
\ni
The first task is to obtain uniform in $\eps$ a priori estimates for the solution of \fref{schrod}, \fref{cauchy}, \fref{poisson}, which are of course crucial in the subsequent nonlinear analysis. Due to the presence of the singular $\frac{1}{\eps^2}$ and $\frac{1}{\eps}$ terms in \fref{schrod}, this task is not obvious here. In subsection \ref{sectpre}, we introduce a well adapted functional framework: a Sobolev scale based on the operators $-\Delta_{x,y}$ and $H_z$. More precisely, for all $m\in \NN$, we introduce the Hilbert space
\be
\label{Bell}
B^m=\left\{u:\,\,\|u\|_{B^m}^2=\|u\|_{L^2(\RR^3)}^2+\|(-\Delta_{x,y})^{m/2}u\|^2_{L^2(\RR^3)}+\|H_z^{m/2}u\|^2_{L^2(\RR^3)}<+\infty\right\}.
\ee
In subsection \ref{sectpre}, we give some equivalent norms which are easier to handle here. Then in subsection \ref{sectapriori} we take advantage of this functional framework and derive some a priori estimates for \fref{schrod}, \fref{cauchy}, \fref{poisson}.

\bs
\ni
{\em Step 2: the filtered system}.

\ms
\ni
In \cite{BMSW1,BACM}, the asymptotics of NLS equations under the form
\be
\label{nls}
i\pa_t u^\eps=\frac{1}{\eps^2}\,H_zu^\eps -\Delta_{x,y} u^\eps +{\mathcal F}(|u^\eps|^2)u^\eps,
\ee
such as the Gross-Pitaevskii equation, was analyzed. In \fref{nls}, ${\mathcal F}:\RR_+\mapsto \RR$ is a given function and the nonlinearity depends locally on the density $|u^\eps|$. It appeared in \cite{BACM} that a fruitful strategy is to filter out the oscillations in time induced by the term $\frac{1}{\eps^2}\,H_z$, without projecting on the eigenmodes of $H_z$\,. Indeed, projecting \fref{nls} on the Hilbert basis $\chi_p$ leads to difficult problems of series summations and of small denominators in oscillating phases. Introducing  the new unknown:
$$v^\eps(t,x,z)=\exp\left(itH_z/\varepsilon^2\right)u^\eps(t,x,z),$$
the filtered system associated to \fref{nls} reads
\be
\label{nlsfiltre}
i\pa_t v^\eps= -\Delta_{x,y} v^\eps +e^{itH_{z}/\varepsilon^2}{\mathcal F}\left(\left|e^{-itH_{z}/\varepsilon^2}v^\eps\right|^2\right)e^{-itH_{z}/\varepsilon^2}v^\eps
\ee
where we used the fact that $H_z$, thus $e^{itH_{z}}$,  commutes with $\pa_x$ and $\pa_y$. Then, the analysis of the limit $\eps\to 0$ amounts to prove that it is possible to define an average of the nonlinearity in \fref{nlsfiltre} with respect to the fast variable $t/\eps^2$.

Let us adapt this strategy to our problem. Introduce
$$\Phi^\eps(t,x,z)=\exp\left(itH_z/\varepsilon^2\right)\Pe(t,x,z).$$
One deduces from \fref{schrod}, \fref{cauchy}, \fref{poisson} the following equation for $\Phe$:
\begin{equation}
\label{Sfiltree}
i\partial_{t}\Phe=-\frac{2B}{\eps}\left(e^{itH_{z}/\varepsilon^2}ze^{-itH_{z}/\varepsilon^2}\right)(i\pa_x\Phe)-\Delta_{x,y}\Phe+ F\left(\frac{t}{\eps^2},\Phi^\eps(t)\right),
\end{equation}
where we introduced the nonlinear function
\be
\label{Pfiltree}
(\tau,u)\mapsto F\left(\tau,u\right)=e^{i\tau H_z}\left(\frac{1}{4\pi r^\eps}\ast \left|e^{-i\tau H_z}u\right|^2\right)e^{-i\tau H_z}u,
\ee
and where $\re$ is still defined by \fref{re}.

\bs
\ni
{\em Step 3: approximation by an intermediate system.}

\ms
\ni
Before performing the limit $\eps\to 0$ in \fref{Sfiltree}, we remark that \fref{Pfiltree} can be approximated in order to get rid of the fast time variable $t/\eps^2$ in the nonlinear term of \fref{Sfiltree}. By writing {\em formally}
\be
\label{formal}
\frac{1}{\sqrt{x^2+y^2+\eps^2z^2}}=\frac{1}{\sqrt{x^2+y^2}}+o(1),
\ee
we remark that
\bee
&\ds\frac{1}{r^\eps}\ast \left|e^{-i\tau H_z}u\right|&=\frac{1}{\sqrt{x^2+y^2}}\ast \left\langle\left|e^{-i\tau H_z}u\right|^2\right\rangle+o(1)\\
&&=\frac{1}{\sqrt{x^2+y^2}}\ast \left\langle\left|u\right|^2\right\rangle+o(1),
\eee
where the symbole $*$ denotes here a convolution in the $(x,y)$ variables only, and where we used the fact that $e^{i\tau H_z}$ is unitary on $L^2_z(\RR)$. Hence, inserting this Ansatz in \fref{Pfiltree} yields
\bee
&F\left(\tau,u\right)&=e^{i\tau H_z}\left(\frac{1}{4\pi\sqrt{x^2+y^2}} \ast \left\langle\left|u\right|^2\right\rangle\right)e^{-i\tau H_z}u+o(1)\\
&&=\left(\frac{1}{4\pi\sqrt{x^2+y^2}}\ast \left\langle\left|u\right|^2\right\rangle\right)u+o(1).
\eee
Denoting
\be
\label{F0}
F_0(u)=\left(\frac{1}{4\pi\sqrt{x^2+y^2}}\ast \left\langle\left|u\right|^2\right\rangle\right)u,
\ee
and introducing the solution $\widetilde \Phe$ of the following intermediate system:
\begin{equation}
\label{interm}
i\partial_{t}\widetilde\Phe=-\frac{2B}{\eps}\left(e^{itH_{z}/\varepsilon^2}ze^{-itH_{z}/\varepsilon^2}\right)(i\pa_x\widetilde\Phe)-\Delta_{x,y}\widetilde\Phe+ F_0\left(\widetilde\Phi^\eps(t)\right),
\end{equation}
we expect that the solution $\Psi^\eps$ of \fref{Sfiltree} satisfies
\be
\label{formal2}
\Phe=\widetilde \Phe +o(1).
\ee
Subsection \ref{sectinterm} is devoted to the rigorous proof of this heuristics. We give sense to the $o(1)$ in Lemma \ref{Poisson} and we prove that the solutions of the two nonlinear equations \fref{Sfiltree} and \fref{interm} are close together and that \fref{formal2} holds true in the sense of the $B^1$ norm. This statement is given in Proposition \ref{prop1}.

\bs
\ni
{\em Step 4: second order averaging of oscillating systems.}

\ms
\ni
Thanks to this Step 3, we can consider the simplest system \fref{interm} instead of \fref{Sfiltree}. We are now left with the analysis of the asymptotics of this intermediate system as $\eps \to 0$. Note that \fref{interm} is under the general form
\begin{equation}
\label{ode}
i\partial_{t}u=\frac{1}{\varepsilon}f\left(\frac{t}{\varepsilon^2}\right)u(t)+g(u(t))
\end{equation}
with
\begin{equation*}
f(\tau)=-2Be^{i\tau H_{z}}ze^{-i\tau H_{z}}i\pa_x\quad \mbox{ and }\quad g(u)=-\Delta_{x,y}u+F_0(u).
\end{equation*}
At this point, a critical fact has to be noticed. Equations under the form
\begin{equation}
\label{ode2}
i\partial_{t}u=f\left(\frac{t}{\varepsilon^2}\right)u(t)+g(u(t))
\end{equation}
can be averaged when, due to some ergodicity property, one can give a sense to the time average
\be
\label{ta}
f^0=\lim_{T\to +\infty}\frac{1}{T}\int_0^Tf(\tau)\,d\tau.
\ee
Indeed, under rather general assumptions, the techniques of averaging of dynamical systems -- see the reference book on the topic by Sanders and Verhulst \cite{SV}-- enable to show that \fref{ode2} is well approximated by the averaged equation
$$
i\partial_{t}u=f^0u(t)+g(u(t)).
$$
Yet, the oscillating term in \fref{ode}, compared to the same term in \fref{ode2}, is multiplied by $\frac{1}{\eps}$. Therefore, a necessary condition in order to perform the averaging of \fref{ode} is that the average $f^0$ of $f$ is zero. In our case, the integral kernel of the operator $e^{i\tau H_{z}}ze^{-i\tau H_{z}}$, defined by
$$\forall u,\qquad e^{i\tau H_{z}}ze^{-i\tau H_{z}}u=\int_\RR G(\tau,z,z')u(z')dz',$$
is given by 
\bee
&G(\tau,z,z')&=\sum_{p\in\NN}\sum_{q\in \NN}e^{i\tau(E_p-E_q)}\left\langle z\chi_p\chi_q\right\rangle\chi_p(z)\chi_q(z')\\
&&=\sum_{p\in\NN}\sum_{q\neq p}e^{i\tau(E_p-E_q)}\left\langle z\chi_p\chi_q\right\rangle\chi_p(z)\chi_q(z').
\eee
In the last inequality, we used the fact that, by Assumption \ref{confinement}, $V_c$ is even. Indeed, this property implies that, for all $p$, $(\chi_p)^2$ is also even, thus $\langle z(\chi_p)^2\rangle=0$. Consequently, since $p\neq q$ implies $E_p\neq E_q$, the kernel $G(\tau,z,z')$ is a series of functions which all have a vanishing average in time. We thus expect that the operator-valued function $f(\tau)$ has the same property:
$$f^0=\lim_{T\to +\infty}\frac{1}{T}\int_0^Tf(\tau)\,d\tau=0.$$
In such a situation,  the theory of averaging has to be pushed to the second order \cite{SV} in order to obtain the limit of \fref{ode} as $\eps\to 0$. Section \ref{sectaver} is devoted to this question of second order averaging, which leads to the limit system \fref{limitS}, \fref{limitP}. The main result of this Section \ref{sectaver} is Proposition \ref{propaver}.

\bs
In the short last Section \ref{lastsect}, we prove our main Theorem \ref{theo} by just gathering the results proved in the previous sections.

\section{The nonlinear analysis}
\label{sect2}

In this section, we obtain some a priori estimates uniform in $\eps$ for the initial system \fref{schrod}, \fref{cauchy}, \fref{poisson} and we prove that it can be approximated by an intermediate system, where we regularize the initial data and where we replace the Poisson nonlinearity by its formal limit given in \fref{F0} . This intermediate system takes the form
\begin{equation}
\label{schrodI}
i\partial_{t}\widetilde\Pe=\frac{1}{\eps^2}H_z\widetilde\Pe-\frac{1}{\eps}2iBz\pa_x\widetilde\Pe-\Delta_{x,y}\widetilde\Pe+W^\eps\widetilde\Pe\,,
\end{equation}
\begin{equation}
\label{cauchyI}
\widetilde\Pe(0,x,y,z)=\widetilde{\Psi_0}(x,y,z),
\end{equation}
\begin{equation}
\label{poissonI}
W^\eps(t,x,z)=\frac{1}{4\pi \sqrt{x^2+y^2}}\ast\left\langle|\widetilde\Pe|^2\right\rangle.
\end{equation}
Notice that \fref{poissonI} is nothing but the Poisson equation \fref{poisson} where we replace $r^\eps=\sqrt{x^2+y^2+\eps^2z^2}$ by $r^0=\sqrt{x^2+y^2}$. Moreover, the initial datum $\widetilde {\Psi_0}$ in \fref{cauchyI} will be chosen as a regularization in $B^m$ of the initial datum $\Psi_0$. Recall the definition \fref{Bell} of the space $B^m$. The main result of this section is the following proposition.
\begin{proposition}[Approximation of the initial system]
\label{prop1}
Assume that $V_c$ satisfies Assumptions \pref{confinement}, \pref{ass2} and that $\Psi_0\in B^1$. For all $\eps\in (0,1]$, denote by $\Pe\in C^0(\RR_+,B^1)$ the unique global solution of the initial system \fref{schrod}, \fref{cauchy}, \fref{poisson}. Then the following holds true.\\
(i) There exists a maximal positive time such that $\Pe$ is bounded uniformly in $\eps$ : the quantity
\be
\label{defT0}
T_0:=\sup\left\{T\geq 0:\,\sup_{\eps\in (0,1]}\|\Psi^\eps\|_{C^0([0,T],B^1)}<+\infty\right\}.
\ee
satisfies $T_0\in (0,+\infty]$. If $T_0<+\infty$ then $$\underset{\footnotesize \eps\to 0}{\mbox{\rm lim sup}}\|\Psi^\eps\|_{C^0([0,T_0],B^1)}=+\infty.$$
(ii) For all $T\in(0,T_0)$, where $T_0$ is defined by \fref{defT0}, for all $\delta>0$ and for all integer $m\geq 2$, there exist $\widetilde{\Psi_0}\in B^m$ and $\eps_\delta$ such that the following holds true.  For all $\eps\in (0,\eps_\delta]$, the intermediate system \fref{schrodI}, \fref{cauchyI}, \fref{poissonI} admits a unique solution $\widetilde\Pe \in C^0([0,T],B^m)$ satisfying the following uniform estimates:
\bea
\label{approx1}
\forall \eps\leq \eps_\delta\qquad \|\Pe -\widetilde\Pe\|_{C^0([0,T],B^1)}&\leq& \delta\\
\label{approx2}
\|\widetilde \Pe\|_{C^0([0,T],B^m)}&\leq& C(\|\Psi_0\|_{B^1})\|\widetilde {\Psi_0}\|_{B^m}\,.
\eea
\end{proposition}
\begin{remark}
It is a priori not excluded that $T_{0}<+\infty$. Indeed, although we are in a repulsive case, the energy conservation does not enable to obtain $\eps$-independant a priori estimates in $B^1$ (see the proof of Lemma \ref{uniformesti}). This may be linked to the possible formation of caustics, as for the nonlinear Schr\"odinger equation in semiclassical regime, see e.g. \cite{Carles}.
\end{remark}

\subsection{Preliminaries}
\label{sectpre}
As we explained in subsection \ref{sketch}, our nonlinear analysis will deeply rely on the use of the functional spaces $B^m$ defined by \fref{Bell} and adapted to the operators $H_z$ and $-\Delta_{x,y}$. The following result was proved in \cite{BACM} by using an appropriate Weyl-H\"ormander pseudodifferential calculus, inspired by \cite{BC,HN}:
\begin{lemma}[\cite{BACM}]
\label{lemBACM}
Under Assumption \pref{confinement}, consider the Hilbert space $B^m$ defined by \fref{Bell} for $m\in \NN$. Then the norm $\|\cdot\|_{B^m}$ in \fref{Bell} is equivalent to the following norm:
\be
\label{Bell2}
\|u\|_{H^m(\RR^3)}+\|V_c(z)^{m/2}u\|_{L^2(\RR^3)}.
\ee
Moreover, for all $u\in B^{m+1}$, we have
\be
\label{derivation}
\|H_z^{1/2} u\|_{B^m}+\|\pa_x u\|_{B^m}+\|\pa_y u\|_{B^m}+\|\pa_z u\|_{B^m}+\|\sqrt{V_c} u\|_{B^m}\lesssim \|u\|_{B^{m+1}}.
\ee
\end{lemma}

The operator $\Delta_{x,y}$ commutes with the rapidly oscillating operator $e^{\pm itH_z/\eps^2}$ and with the operator  $iz\pa_x$. This will enable us to obtain uniform bounds for the solution of \fref{schrod} by simply applying $\Delta_{x,y}$ to this equation. Unfortunately, the operator $H_z$ does not satisfy this property. For this reason, we introduce the following operator:
\be
\label{HH}
H_\eps=H_z-2i\eps Bz\pa_x-\eps^2\pa^2_x=-\pa_z^2+V_c(z)+(i\eps\pa_x-Bz)^2\,.
\ee
This operator enables to define another norm equivalent to the $B^m$ norm. The following lemma is proved in the Appendix \ref{appZ}.
\begin{lemma}
\label{sobo2}
The operator $H_\eps$ defined by \fref{HH} on $L^2(\RR^3)$ with domain $B^2$ is self-adjoint and nonnegative. There exists a constant $C_1>0$ such that, for all $\eps\in (0,1]$ and for all $u\in B^1$, we have
\be
\label{Bell31}
\frac{1}{C_1}\|u\|_{B^1}^2\leq \|u\|_{L^2(\RR^3)}^2+\|(-\Delta_{x,y})^{1/2}u\|^2_{L^2(\RR^3)}+\|H_\eps^{1/2}u\|^2_{L^2(\RR^3)}\leq C_1\|u\|_{B^1}^2\,.
\ee
Moreover, for all integer $m\geq 2$, there exists $\eps_m\in (0,1]$ such that, for all $\eps\in (0,\eps_m]$, for all $u\in B^m$, we have
\be
\label{Bell3}
\frac{1}{2}\|u\|_{B^m}^2\leq \|u\|_{L^2(\RR^3)}^2+\|(-\Delta_{x,y})^{m/2}u\|^2_{L^2(\RR^3)}+\|H_\eps^{m/2}u\|^2_{L^2(\RR^3)}\leq 2\|u\|_{B^m}^2\,.
\ee
\end{lemma}

\subsection{A priori estimates}
\label{sectapriori}
In this subsection, we obtain a priori estimate uniform in $\eps$ for the initial Schr\"odinger-Poisson model \fref{schrod}, \fref{poisson} and the intermediate model \fref{schrodI}, \fref{cauchyI}, \fref{poissonI}. Remark first that these two models can be considered in a unified way. For all $u\in B^1$ and for $\alpha\in \{0,1\}$, denote
\be
\label{unif}
F_\alpha(u)=\left(\frac{1}{4\pi \sqrt{x^2+y^2+\alpha \eps^2z^2}}\ast \left(|u|^2\right)\right)u\,,
\ee
where the convolution holds on the three variables $(x,y,z)\in \RR^3$. Remark that for $\alpha=0$, this definition coincides with the definition \fref{F0}. We shall consider for $\eps\in (0,1]$ and $\alpha\in\{0,1\}$ the nonlinear equation
\begin{equation}
\label{schrodgene}
i\partial_{t}u^\eps=\frac{1}{\eps^2}H_\eps u^\eps-\pa^2_yu^\eps+F_\alpha(u^\eps)\,,
\end{equation}
\begin{equation}
\label{cauchygene}
u^\eps(0,x,y,z)=u_0(x,y,z),
\end{equation}
where the operator $H_\eps$ was defined by \fref{HH}. Note that for $u_0=\Psi_0$ and $\alpha=1$, \fref{schrodgene}, \fref{cauchygene} is the initial system \fref{schrod}, \fref{cauchy}, \fref{poisson}, and that for $u_0=\widetilde{\Psi_0}$ and $\alpha=0$, \fref{schrodgene}, \fref{cauchygene} is the intermediate system \fref{schrodI}, \fref{cauchyI}, \fref{poissonI}. Let us first state a technical lemma concerning the nonlinearities $F_1$ and $F_0$, which is proved in Appendix \ref{appA}.
\begin{lemma}
\label{tame estimate}
There exists a constant $C>0$ such that, for all $\eps\in (0,1]$, for $\alpha=0$ or 1, we have
\be
\label{estB1}
\forall u,v\in B^1,\quad\left\|F_\alpha(u)-F_\alpha(v)\right\|_{B^1}\leq C\left(\|u\|_{B^1}^2+\|v\|_{B^1}^2\right)\|u-v\|_{B^1}\,,
\ee
where $F_\alpha$ is defined by \fref{unif}.
Moreover, for all $m\in \NN^*$, there exists $C_m>0$ such that we have the tame estimate
\be
\label{tame}
\forall \eps \in (0,1],\,\forall \alpha\in\{0,1\},\,\forall u\in B^m,\quad \left\|F_\alpha(u)\right\|_{B^m}\leq C_m\|u\|_{B^1}^2\,\|u\|_{B^m}\,.
\ee
\end{lemma}
Now we are able to derive uniform a priori estimates for the solution of \fref{schrodgene}, \fref{cauchygene}.
\begin{lemma}
\label{uniformesti}
Let $\eps\in (0,1]$, $\alpha\in\{0,1\}$ and $u_0\in B^1$. Then the solution $u^\eps$ of the equation \fref{schrodgene}, \fref{cauchygene} exists and is unique in $C^0([0,+\infty),B^1)$ and the following uniform in $\eps$ estimates hold true.\\
(i) For all $M>0$, there exist $T>0$, only depending on $M$ and $\|u_0\|_{B^1}$, such that, for all $\eps\in (0,1]$, we have
\be
\label{estimaB1}
\|u^\eps\|_{C^0([0,T],B^1)}\leq (1+M)\|u_0\|_{B^1}\,.
\ee
(ii) Let $m\geq 2$ an integer and assume that $u_0\in B^m$. Then, for all $\widetilde T>0$, we have the estimate
\be
\label{estimaBm}
\forall \eps\in (0,\eps_m],\qquad \|u^\eps\|_{C^0([0,\widetilde T],B^m)}\leq C\|u_0\|_{B^m}\exp\left(C\widetilde T\|u^\eps\|_{C^0([0,T],B^1)}^2\right).
\ee
where $\eps_m>0$ is as in Lemma \pref{sobo2}.
\end{lemma}
\begin{proof}

\bs
\ni
{\em Step 1: the Cauchy problem and the conservation laws.} For any given $\eps>0$, the existence and uniqueness of a maximal solution $u^\eps\in C^0([0,\overline{T}),B^1)$ can be obtained by standard techniques \cite{CAZ}.  We leave this first part of the proof to the reader. This solution satisfies both $L^2$ and energy conservation laws:
\be
\label{consgene}
\forall t\geq 0,\quad \|u^\eps(t)\|_{L^2}=\|u_0\|_{L^2}\quad\mbox{and}\quad \calE_\alpha(u^\eps(t))=\calE_\alpha(u_0),
\ee
where the energy $\calE_\alpha$ is defined by
\bee
&&\calE_\alpha(u)=\frac{1}{\eps^2}(H_\eps u,u)_{L^2}+\|\pa_y u\|^2_{L^2}+\frac{1}{2}(F_\alpha(u),u)_{L^2}\\
&&=\frac{1}{\eps^2}\|\pa_z u\|_{L^2}^2+\frac{1}{\eps^2}\|\sqrt{V_c}u\|_{L^2}^2+\frac{1}{\eps^2}\|(\eps\pa_x+iBz)u\|_{L^2}^2+\|\pa_yu\|_{L^2}^2+\frac{1}{2}(F_\alpha(u),u)_{L^2}\,.
\eee
We recall that the operator $H_\eps$ is defined by \fref{HH}. These conservation laws show that the solution $u^\eps$ is global, ie that $\overline T=+\infty$. Unfortunately, due to the $\frac{1}{\eps^2}$ terms in this expression, one cannot use the energy conservation to get uniform in $\eps$ estimates. Instead, we will directly write the equations satisfied by $\partial_{x}u^{\eps}, \partial_{y}u^{\eps}$ or $(H_{\eps})^{1/2}u^{\eps}$ and use the standard $L^2$-estimates for these equations and the fact that the self-adjoint operators $H_\eps$, $\partial_{x}$ and $\partial_{y}$ commute together.

\bs
\ni
{\em Step 2: $B^1$ estimate.} This yields
$$
i\partial_{t}(\na_{x,y}u^\eps)(t)=\frac{1}{\eps^2}H_{\eps}(\nabla_{x,y}u^{\eps})-\partial_{y}^2(\nabla_{x,y}u^{\eps})+\na_{x,y}\left(F_\alpha(u^\eps)\right)
$$
and
$$
i\partial_{t}\left(H_\eps^{1/2}u^\eps\right)(t)=\frac{1}{\eps^2}H_{\eps}(H_{\eps}^{1/2}u^{\eps})-\partial_{y}^2(H_{\eps}^{1/2}u^{\eps})+H_\eps^{1/2}\left(F_\alpha(u^\eps)\right).
$$
Hence, 
$$\begin{array}{r}
\ds\|u^\eps(t)\|_{L^2}+\|\na_{x,y} u^\eps(t)\|_{L^2}+\|H_\eps^{1/2}u^\eps(t)\|_{L^2}\leq \|u_0\|_{L^2}+\|\na_{x,y} u_0\|_{L^2}+\|H_\eps^{1/2}u_0\|_{L^2}\\[3mm]
\ds+C\int_0^t\left(\|\na_{x,y} F_\alpha (u^\eps(s))\|_{L^2}+\|H_\eps^{1/2}F_\alpha(u^\eps(s))\|_{L^2}\right)\,ds
\end{array}
$$
and, for $\eps\in (0,1]$, the equivalence of norms given in Lemma \ref{sobo2}, yields
\bea
&\|u^\eps(t)\|_{B^1}&\leq C\|u_0\|_{B^1}+C\int_0^t \|F_\alpha (u^\eps(s))\|_{B^1}\,ds\nonumber\\
&&\leq C\|u_0\|_{B^1}+C\int_0^t \|u^\eps(s)\|^3_{B^1}\,ds,\label{gron2}
\eea
where we used \fref{estB1} with $v=0$ to estimate $F_\alpha (u^\eps(s))$.  Hence, by applying the Gronwall lemma to the integral inequality \fref{gron2}, we prove Item {\em (i)} of the Lemma.

\bs
\ni
{\em Step 3: $B^m$ estimate.} Let $T>0$, $m\geq 2$, $u_0\in B^m$ and let $\eps\in (0,\eps_m]$, where $0<\eps_m\leq 1$ as in Lemma \ref{sobo2}. Since the operators $H_\eps$ and $\Delta_{x,y}$ commute together, $H_{\eps}^{m/2}u^{\eps}$ satifies the following equation:
$$
i\partial_{t}\left(H^{m/2}_\eps u^\eps\right)(t)=\frac{1}{\eps^2}H_{\eps}(H_{\eps}^{m/2}u^{\eps})-\partial_{y}^2(H_{\eps}^{m/2}u^{\eps})+H_\eps^{m/2} \left(F_\alpha(u^\eps)\right),
$$
thus, for all $t\in [0,T]$,
\bea
&\|H^{m/2}_\eps u^\eps(t)\|_{L^2}&\leq \|H^{m/2}_\eps u_0\|_{L^2}+\int_0^t\|H^{m/2}_\eps \left(F_\alpha(u^\eps(s))\right)\|_{L^2}\,ds,\nonumber\\
&&\leq C\|u_0\|_{B^m}+C\int_0^t \|F_\alpha(u^\eps(s))\|_{B^m}\,ds\nonumber\\
&&\leq  C\|u_0\|_{B^m}+C\|u^\eps\|_{C^0([0,T],B^1)}^2\int_0^t \|u^\eps(s)\|_{B^m}\,ds,\label{bor1}
\eea
where we used Lemma \ref{sobo2} and the tame estimate \fref{tame}. Similarly, $-\Delta_{x,y}u^{\eps}$ satisfies the following equation:
$$i\partial_{t}(-\Delta_{x,y}u^{\eps})(t)=\frac{1}{\eps^2}H_{\eps}(-\Delta_{x,y}u^{\eps})-\partial_{y}^2(-\Delta_{x,y}u^{\eps})-\Delta_{x,y}\left(F_{\alpha}(u^{\eps)}\right)$$
and, using the definition of $B^{m}$ \fref{Bell} and \fref{tame} yields:
\bea
&\|(-\Delta_{x,y})^{m/2} u^\eps(t)\|_{L^2}&\leq \|(-\Delta_{x,y})^{m/2} u_0\|_{L^2}+\int_0^t\|(-\Delta_{x,y})^{m/2} \left(F_\alpha(u^\eps(s))\right)\|_{L^2}\,ds,\nonumber\\
&&\leq  C\|u_0\|_{B^m}+C\|u^\eps\|_{C^0([0,T],B^1)}^2\int_0^t \|u^\eps(s)\|_{B^m}\,ds.\label{bor2}
\eea
Therefore, by using again the equivalence of norms given by Lemma \ref{sobo2} and the $L^2$ conservation law in \fref{consgene}, we deduce from \fref{bor1} and \fref{bor2} that, for $t\leq T$, we have
$$\|u^\eps(t)\|_{B^m}\leq  C\|u_0\|_{B^m}+C\|u^\eps\|_{C^0([0,T],B^1)}^2\int_0^t \|u^\eps(s)\|_{B^m}\,ds,$$
and the Gronwall lemma gives \fref{estimaBm}.
\end{proof}

\subsection{Proof of Proposition \ref{prop1}}
\label{sectinterm}
In this subsection, we prove Proposition \ref{prop1}, ie we show that this solution can be uniformly approximated by a regular solution of the intermediate system. We first state a technical lemma on the Poisson kernels, which is proved in the Appendix \ref{appB}.
\begin{lemma}
\label{Poisson}
There exists a constant $C>0$ such that, for all $\eps \in(0,1]$, we have
\be
\label{diffnoyau}
\forall u\in B^2,\qquad \left\|F_1(u)-F_0(u)\right\|_{B^1}\leq C\,\eps^{1/3}\,\|u\|_{B^2}^3\,,
\ee
where $F_0$ and $F_1$ are defined by \fref{unif}.
\end{lemma}

\ni
We are now ready to prove the main result of this section.

\bs
\ni
{\em Proof of Proposition \pref{prop1}.} Let $\Psi_0\in B^1$, let an integer $m\geq 2$ be fixed, and define the regularized initial datum $\widetilde {\Psi_0}$ by
\be
\label{reginit}
\widetilde {\Psi_0}=\left(I-\eta\Delta_{x,y}\right)^{-m/2}\left(I+\eta H_z\right)^{-m/2}\Psi_0\,,
\ee
where $\eta>0$ is a small parameter that will be fixed further and where $I$ denotes the identity operator. Denote by $\Pe$ the solution of the initial system \fref{schrod}, \fref{cauchy}, \fref{poisson} and by $\widetilde\Pe$ the solution of the intermediate system \fref{schrodI}, \fref{cauchyI}, \fref{poissonI} with the initial datum \fref{reginit}. We shall estimate the difference $\Pe-\widetilde \Pe$.

\ms
\ni
{\em Step 1: uniform bounds for $\Pe$}. Let $0<\eps \leq 1$. From  Lemma \ref{uniformesti} {\em (i)}, we first deduce that there exists $T_1>0$ only depending on $\|\Psi_0\|_{B^1}$ such that, for all $\eps\in(0,1]$
$$
\|\Pe\|_{C^0([0,T_1],B^1)}\leq 2\|\Psi_0\|_{B^1}.
$$
This implies that $T_0$ defined by \fref{defT0} satisfies $T_0\geq T_1>0$. Clearly, if $T_0<+\infty$, we have
$$\underset{\footnotesize \eps\to 0}{\mbox{\rm lim sup}}\|\Psi^\eps\|_{C^0([0,T_0],B^1)}=+\infty,$$
otherwise by reiterating the above procedure we could find a uniform bound on $[0,T_2]$ with $T_2>T_0$.

Now we fix $T\in(0,T_0)$ and $\delta>0$ for the sequel of this proof. Definition \fref{defT0} of $T_0$ implies that
\be
\label{estimtousB1bis}
\|\Pe\|_{C^0([0,T],B^1)}\leq C\left(\|\Psi_0\|_{B^1}\right), \mbox{ independent of }\eps\in (0,1].
\ee

\ms
\ni
{\em Step 2: bounds for the initial datum $\widetilde {\Psi_0}$}. First, we deduce from \fref{reginit} that
$$(I-\Delta_{x,y}+H_z)^{1/2}\widetilde {\Psi_0}=(I-\eta\Delta_{x,y})^{-m/2}(I+\eta H_z)^{-m/2}(I-\Delta_{x,y}+H_z)^{1/2}{\Psi_0}\,,$$
hence
\bee
&&\|(I-\Delta_{x,y}+ H_z)^{1/2}\widetilde {\Psi_0}\|_{L^2}\\
&&\qquad \leq \|(I-\eta\Delta_{x,y})^{-m/2}(I+\eta H_z)^{-m/2}(I-\Delta_{x,y}+ H_z)^{1/2}{\Psi_0}\|_{L^2}\\
&&\qquad \leq \|(I-\Delta_{x,y}+ H_z)^{1/2}{\Psi_0}\|_{L^2}
\eee
where we used the fact that the operators $(I-\eta\Delta_{x,y})^{-m/2}$ and $(I+\eta H_z)^{-m/2}$ are bounded on $L^2$, with  bounds equal to 1. Therefore, using \fref{equivB11}, we obtain
\be
\label{estiminit0}
\|\widetilde {\Psi_0}\|_{B^1}\leq \|\Psi_0\|_{B^1},
\ee
where we recall that the right-hand side is independent of $\eps$. 

Next, we get from \fref{reginit} the two following identities: for all integer $\ell\leq m,$
$$(-\Delta_{x,y})^{\ell/2+1/2}\widetilde {\Psi_0}=(-\Delta_{x,y})^{\ell/2}(I-\eta\Delta_{x,y})^{-\ell/2}(I-\eta\Delta_{x,y})^{\ell/2-m/2}(I+\eta H_z)^{-m/2}(-\Delta_{x,y})^{1/2}{\Psi_0}\,,$$
and
$$H_z^{\ell/2+1/2}\widetilde {\Psi_0}=H_z^{\ell/2}(I+\eta H_z)^{-\ell/2}(I+\eta H_z)^{\ell/2-m/2}(I-\eta \Delta_{x,y})^{-m/2}H_z^{1/2}{\Psi_0}\,.$$
Thus, from the bound
$$\forall \lambda\in\RR_+,\quad \lambda^{\ell/2}(1+\eta\lambda)^{-{\ell/2}}\leq C\eta^{-\ell/2}\,,$$
we deduce that both operators $(-\Delta_{x,y})^{\ell/2}(I-\eta\Delta_{x,y})^{-\ell/2}$ and $H_z^{\ell/2}(I+\eta H_z)^{-\ell/2}$ are bounded on $L^2$, with bounds equal to $C\eta^{-\ell/2}$, and thus 
\be
\label{estiminit}
\forall \ell\leq m,\qquad \|\widetilde {\Psi_0}\|_{B^{\ell+1}}\leq C\eta^{-\ell/2}\,\|\Psi_0\|_{B^1},
\ee
where we recall the definition \fref{Bell} of the $B^m$ norms.

Finally, we obtain also from \fref{reginit} that
$$(I-\Delta_{x,y}+H_z)^{1/2}(\Psi_0-\widetilde {\Psi_0})=\left(I-(I-\eta\Delta_{x,y})^{-m/2}(I+\eta H_z)^{-m/2}\right)(I-\Delta_{x,y}+H_z)^{1/2}{\Psi_0}\,.$$
Decompose $v=(I-\Delta_{x,y}+H_z)^{1/2}{\Psi_0}$ on the Hilbert basis $(\chi_p)_{p\in\NN}$ of eigenmodes of $H_z$:
$$v(x,y,z)=\sum_{p\in\NN}v_p(x,y)\,\chi_p(z)$$
and denote by $\widehat{v_p}(\xi)$, $\xi\in\RR^2$, the Fourier transform of $v_p(x,y)$. By \fref{equivB11}, we have
$$\|\Psi_0-\widetilde {\Psi_0}\|^2_{B^1}=\sum_{p\in\NN}\int_{\RR^2}\left(1-(1+\eta |\xi|^2)^{-m/2}(1+\eta E_p)^{-m/2}\right)^2|\widehat{v_p}(\xi)|^2\,d\xi\,.$$
Hence, using that
\be
\label{four}
\sum_{p\in\NN}\int_{\RR^2}|\widehat{v_p}(\xi)|^2\,d\xi=\|\Psi_0\|_{B^1}^2<+\infty
\ee
and that
$$\forall \xi\in \RR^2,\,\forall p\in \NN,\qquad \lim_{\eta\to 0}\left(1-(1+\eta |\xi|^2)^{-m/2}(1+\eta E_p)^{-m/2}\right)=0,$$
we deduce from Lebesgue's dominated convergence theorem and from the convergence of the series in \fref{four} that 
\be
\label{convinit}
\lim_{\eta\to 0}\|\Psi_0-\widetilde {\Psi}_0\|_{B^1}=0.
\ee

\ms
\ni
{\em Step 3: uniform a priori estimates for $\widetilde \Pe$}. Consider
\be
\label{estimtousB1}
T_{\eta}:=\sup\{\tau\in (0,T]\,:\ \forall \eps\in (0,1],\ \|\widetilde \Pe\|_{C^0([0,T_\eta],B^1)}\leq 2\|\Pe\|_{C^0([0,T],B^1)}\}.
\ee
Note that, from \fref{estiminit0} and Lemma \ref{uniformesti} {\em (i)}, we know that $T_\eta\in (0,T]$ is well-defined. Then, from Lemma \ref{uniformesti}  {\em (ii)}, we deduce the following estimate:
\bea
&\forall \eps\in (0,\eps_m],\quad\forall \ell\leq m,\qquad \|\widetilde \Pe\|_{C^0([0,T_\eta],B^{\ell+1})}&\leq C\left(\|\widetilde \Pe\|_{C^0([0,T_\eta],B^1)}\right)\|\widetilde {\Psi_0}\|_{B^{\ell+1}}\nonumber\\
&&\leq C\left(\|\Pe\|_{C^0([0,T],B^1)}\right)\|\widetilde {\Psi_0}\|_{B^{\ell+1}}\nonumber\\
&&\leq C\left(\|\Psi_0\|_{B^1}\right)\|\widetilde {\Psi_0}\|_{B^{\ell+1}}\label{Bmbis}
\eea
where we used \fref{estimtousB1} and \fref{estimtousB1bis}.

\bs
\ni
{\em Step 4: estimate of the difference $\Psi^\eps-\widetilde \Psi^\eps$}. Using the notations defined in \fref{HH} and  \fref{unif}, $\Pe$ and $\widetilde \Pe$ satisfy \fref{schrodgene},\fref{cauchygene} with $\alpha=1$, $u_{0}=\Psi_{0}$ and $\alpha=0$, $u_{0}=\widetilde{\Psi_{0}}$ respectively. The Duhamel formulation of these equations read respectively
$$
\Pe(t)=e^{-it (H_{\eps}-\partial_{y}^2)}\Psi_0+\int_0^te^{-i(t-s)(H_{\eps}-\partial_{y}^2)}F_1(\Pe(s))\,ds,
$$
$$
\widetilde \Pe(t)=e^{-it(H_{\eps}-\partial_{y}^2)}\widetilde{\Psi_0}+\int_0^te^{-i(t-s)(H_{\eps}-\partial_{y}^2)}F_0(\widetilde\Pe(s))\,ds.
$$
Hence, for all $t\in [0,T_\eta]$ and $\eps\in (0,\eps_m]$,
\bee
&\|\Pe(t)-\widetilde\Pe(t)\|_{B^1}&\leq \|\Psi_0-\widetilde{\Psi_0}\|_{B^1}+\int_0^t\|F_1(\Pe(s))-F_1(\widetilde\Pe(s))\|_{B^1}ds\\
&&\quad +\int_0^t\|F_1(\widetilde\Pe(s))-F_0(\widetilde\Pe(s))\|_{B^1}ds\\
&&\leq  \|\Psi_0-\widetilde{\Psi_0}\|_{B^1}+C\int_0^t\|\Pe(s)-\widetilde\Pe(s)\|_{B^1}ds+C\,\eps^{1/3}\eta^{-3/2},\\
\eee
where we used \fref{estB1}, \fref{estimtousB1bis}, \fref{estimtousB1}, \fref{diffnoyau} and \fref{Bmbis} with $\ell=1$, coupled to \fref{estiminit}. Here $C$ denotes a generic constant depending only on $T$ and $\|\Psi_0\|_{B^1}$. Hence, by the Gronwall lemma, we get, for all $t\in [0,T_\eta]$,
\be
\label{gron}
\|\Pe(t)-\widetilde\Pe(t)\|_{B^1}\leq \left( \|\Psi_0-\widetilde{\Psi_0}\|_{B^1}+C\,\eps^{1/3}\eta^{-3/2}\right)e^{CT}.
\ee
Now, according to \fref{convinit}, we fix $\eta$ such that
$$ \|\Psi_0-\widetilde{\Psi_0}\|_{B^1}e^{CT}\leq \min\left(\frac{\delta}{2},\frac{1}{3}\|\Pe\|_{C^0([0,T],B^1)}\right)$$
and, in a second step, we fix $\eps_\delta\in (0,\eps_m]$ such that
$$C\,\eps_\delta^{1/3}\eta^{-3/2}e^{CT}\leq \min\left(\frac{\delta}{2},\frac{1}{3}\|\Pe\|_{C^0([0,T],B^1)}\right).$$
From \fref{gron}, we deduce that
\be
\label{Tdelta1}
\forall t\in [0,T_\eta],\quad \forall \eps\in (0,\eps_\delta],\quad \|\Pe(t)-\widetilde\Pe(t)\|_{B^1}\leq \min\left(\delta,\frac{2}{3}\|\Pe\|_{C^0([0,T],B^1)}\right).
\ee
Therefore, we have
\begin{align}
\nonumber
\|\widetilde \Pe\|_{C^0([0,T_\eta],B^1)}&\leq \| \Pe\|_{C^0([0,T_\eta],B^1)}+\|\Pe-\widetilde \Pe\|_{C^0([0,T_\eta],B^1)}\\
\label{poupoupidou}
&\leq \frac{5}{3}\|\Pe\|_{C^0([0,T],B^1)}.
\end{align}
We claim that $T_{\eta}=T$. Indeed, if $T_{\eta}<T$, then, applying again Lemma \ref{uniformesti} at $T_{\eta}$ and using \fref{poupoupidou} enables to find $\tau>0$ such that, for all $\eps \in (0,1)$,
$$\|\widetilde{\Psi^{\eps}}\|_{C^0([T_{\eta},T_{\eta}+\tau],B^1)}\leq 2 \|\Psi^{\eps}\|_{C^0([0,T],B^1)},$$
which, together with \fref{poupoupidou}, contradicts the definition \fref{estimtousB1} of $T_{\eta}$. Finally, \fref{Tdelta1} gives \fref{approx1} and \fref{Bmbis} with $\ell=m-1$ gives \fref{approx2}. The proof of Proposition \ref{prop1} is complete.
\qed

\section{Second order averaging}
\label{sectaver}

In this section, we focus on the intermediate system \fref{schrodI}, \fref{cauchyI}, \fref{poissonI} as $\eps$ goes to zero. As we explained in subsection \ref{sketch}, it is interesting to consider the filtered version of this equation. Let $\widetilde{\Psi_0}\in B^m$ be a given initial data, let $\widetilde\Pe$ be the corresponding solution of \fref{schrodI}, \fref{cauchyI}, \fref{poissonI} and set
\be
\label{changeinc}
\widetilde {\Phi^\eps}(t,\cdot)=\exp\left(itH_z/\eps^2\right)\widetilde\Pe(t,\cdot).
\ee
This function satisfies the system
\begin{align}
\label{interm2}
&i\partial_{t}\widetilde\Phe=-\frac{2B}{\eps}\left(e^{itH_{z}/\varepsilon^2}ze^{-itH_{z}/\varepsilon^2}\right)(i\pa_x\widetilde\Phe)-\Delta_{x,y}\widetilde\Phe+ F_0\left(\widetilde\Phi^\eps(t)\right),\\ 
\nonumber &\widetilde\Phe(t=0)=\widetilde{\Psi_0},
\end{align}
where $F_0$ is defined by \fref{F0}. The advantage of this intermediate system, compared to \fref{Sfiltree} is that the nonlinearity $F_0(\widetilde\Phe)$ has no dependence in the fast variable $\frac{t}{\eps^2}$.

We will analyze the filtered system \fref{interm2} in the framework of second order averaging of fast oscillating ODEs under the form \fref{ode} --see \cite{SV}--, that we adapt here to our context of nonlinear PDEs.  Recall that $(E_p)_{p\in \NN}$, $(\chi_p)_{p\in \NN}$ are the complete families of eigenvalues and eigenfunctions of the operator $H_z$ and denote by $\Pi_p$ the spectral projector on $\chi_p$:
$$\forall \Phi\in L^2(\RR^3),\qquad \Pi_p\Phi=\langle \Phi \chi_p\rangle \chi_p.$$
Introduce now the following unbounded operator on $L^2(\RR^3)$:
\be
\label{A0}
A_0=-\pa_x^2\sum_{p\geq 0} \alpha_p \,\Pi_p \quad \mbox{with}\quad \alpha_p=1-\sum_{q\neq p}\frac{\left\langle 2Bz \chi_p\chi_q\right\rangle^2}{E_q-E_p}.
\ee
With this notation, the limit system \fref{limitS}, \fref{limitP}, \fref{psilimit} can be rewritten in a more compact form as
\be
\label{eqA0}
i\pa_t\Phi=A_0\Phi-\pa_y^2\Phi+F_0(\Phi),\qquad \Psi(t=0)=\Psi_0.
\ee
We state the main results of this section in the following two propositions.
\begin{proposition}
\label{lemA0}
Assume that $V_c$ satisfies Assumptions \pref{confinement} and \pref{ass2}. Then the following properties hold true.\\
(i) The unbounded operator $A_0$ defined by \fref{A0} on $L^2(\RR^3)$ with the domain
$$D(A_0)=\{\Phi\in L^2(\RR^3):\quad \pa_x^2\sum_{p\geq 0} \alpha_p \,\Pi_p\Phi\in L^2(\RR^3)\}$$
is selfadjoint. Moreover, the operator $A_{0}$ satisfies  
\be
\label{embed}
\forall \ell\geq 0,\ \forall u\in B^{2n_{0}+4+\ell}\qquad \|A_{0}u\|_{B^{\ell}}\leq C \|u\|_{B^{2n_{0}+4+\ell}}
\ee
where $n_0$ is as in Assumption \pref{ass2}.\\
(ii) Let $\Psi_0\in B^1$. The limit system \fref{eqA0} admits a unique maximal solution $\Phi\in C^0([0,T_{max}),B^1)$. If $T_{max}<+\infty$ then $\|\Phi(t)\|_{B^1}\to +\infty$ as $t\to T_{max}$. 
\end{proposition}
\begin{proposition}[Averaging of the intermediate system]
\label{propaver}
Assume that $V_c$ satisfies Assumptions \pref{confinement} and \pref{ass2}. Then there exists an integer $m\geq 2$ such that the following holds true. For $\widetilde {\Psi_0}\in B^m$, we consider the solution $\widetilde {\Phe}\in C^0([0,+\infty),B^m)$  of \fref{interm2} and the maximal solution $\widetilde \Phi\in C^0([0,T_{max}),B^1)$ of the limit system with $\widetilde {\Psi_0}$ as initial data:
\be
\label{lim}
\widetilde \Phi(t)=e^{-it(A_0-\pa^2_y)}\,\widetilde {\Psi_0}-i\int_0^t e^{-i(t-s)(A_0-\pa^2_y)}F_0(\widetilde \Phi(s))ds.
\ee
We assume that there exist $T\in (0,T_{max})$, $\eps_0>0$ such that
\be
\label{estiaver1}
M:=\sup_{\eps\in(0,\eps_0]}\|\widetilde \Phe\|_{C^0([0,T],B^m)}<+\infty.
\ee
Then we have
\be
\label{estili}
\|\widetilde \Phe-\widetilde \Phi\|_{C^0([0,T],B^1)}\leq \eps\,C_M,
\ee
where $C_M$ is independent of $\eps$.
\end{proposition}

\bs
\ni 
\subsection{Well-posedness of the limit system}
In this section, we prove Proposition \ref{lemA0}.

\bs
\ni
{\em Step 1. Basic properties of the operator $A_0$.} First, from $V_c(z)\geq a^2z^2$, we deduce that the $p$th eigenvalue of $H_z$ is larger than the $p$th eigenvalue of the harmonic oscillator $-\frac{d^2}{dz^2}+(a^2+B^2)z^2$:
\be
\label{weyl}
\forall p\in \NN,\qquad E_p\geq \sqrt{a^2+B^2}(2p+1).
\ee
From Assumption \ref{ass2}, we deduce that the coefficients $\alpha_p$ in \fref{A0} satisfy
\bee
&|\alpha_p|&\leq 1+C(1+p)^{n_0}\sum_{q\geq 0}\left\langle 2Bz \chi_p\chi_q\right\rangle^2=1+C(1+p)^{n_0}\|Bz\chi_p\|_{L^2}^2\\
&&\leq CE_p^{n_0+1},
\eee
where we used \fref{weyl} and that $\|Bz\chi_p\|_{L^2}\leq E_p^{1/2}$. Now, consider a nonnegative integer $\ell$ and $u$ in $B^{2n_{0}+4+\ell}$. Let $n_{0}$ be defined as in Assumption \ref{ass2}, and decompose $u$ over the $\chi_{p}$ family which is orthogonal in $L^2$.
\bee
&\|A_{0}u\|^2_{B^\ell}&=\sum_{p\geq 0}\alpha_{p}^2\|\partial_{x}^2\Pi_{p}u\|^2_{B^\ell}\\
&&\leq C\sum_{p\geq 0}E_{p}^{2n_{0}+2}\|\Pi_{p}u\|^2_{B^{\ell+2}}\leq C\sum_{p\geq 0}\|H_{z}^{n_{0}+1}\Pi_{p}u\|^2_{B^{\ell+2}}\\
&&\leq C \sum_{p\geq 0}\|\Pi_{p}u\|^2_{B^{2n_{0}+4+\ell}}=C\|u\|^2_{B^{2n_{0}+4+\ell}}
\eee
where we used Lemma \ref{lemBACM}. This proves \fref{embed}.

Furthermore, by passing to the limit as $N\to +\infty$ in the identity
$$\forall \Phi,\Psi\in D(A),\qquad \sum_{p=0}^N\alpha_p(\pa_x^2\Pi_p\Phi,\Pi_p\Psi)_{L^2}= \sum_{p=0}^N\alpha_p(\Pi_p\Phi,\pa_x^2\Pi_p\Psi)_{L^2},$$
we obtain that the operator $A_0$ is symmetric. Moreover, the equation $A_0\Phi+i\Phi=f$ admits a solution $\Phi\in D(A_0)$ for all $f\in L^2(\RR^3)$. Indeed, the projection of this equation on $\chi_p$ reads
$$-\alpha_p\pa^2_x\phi_p+i\phi_p=f_p$$
and this elliptic equation can obviously be solved for all $f_p\in L^2(\RR^2)$.
Therefore, by the standard criterion for selfadjointness \cite{RS}, the operator $A_0$ is selfadjoint. We have proved the first part of Proposition \ref{lemA0}.

\bs
\ni 
{\em Step 2. Well-posedness  and stability of the limit system.} The operator $A_0$ being selfadjoint, the Stone theorem can be applied and the operator $-iA_0$ generates a unitary group of continuous operators $e^{-iA_0t}$ on $L^2$ and also on $B^1$. The Duhamel formulation of \fref{eqA0} reads
\be
\label{duhalim}
\Phi(t)=e^{-it(A_0-\pa^2_y)}\Psi_0-i\int_0^te^{-i(t-s)(A_0-\pa^2_y)}F_0(\Phi(s))ds
\ee
(recall that $A_0$ and $\pa^2_y$ commute together). Since, by \fref{estB1}, the application $F_0$ is locally Lipschitz continuous on $B^1$, it is easy to prove by a standard fixed point technique that \fref{duhalim} admits a unique maximal solution $\Phi\in C^0([0,T_{max}),B^1)$. The details are left to the reader. Note that, if $T_{max}<+\infty$, then  $\|\Phi(t)\|_{B^1}\to +\infty$ as $t\to T_{max}$. Item {\em (ii)} of Proposition \ref{lemA0} is proved.
\qed

\begin{remark}
\label{remstab}
In fact, this strategy of proof by a fixed point mapping leads to a stability result. For all $\eta>0$ and for all $T\in (0,T_{max})$, there exists $\delta_{\eta,T}>0$ such that the following holds true. For all $\widetilde \Psi_0$ satisfying 
$$\|\Psi_0-\widetilde \Psi_0\|_{B^1}\leq \delta_{\eta,T},$$
 the equation \fref{lim}
\begin{equation*}
\label{duhastab}
\widetilde \Phi(t)=e^{-it(A_0-\pa^2_y)}\widetilde \Psi_0-i\int_0^te^{-i(t-s)(A_0-\pa^2_y)}F_0(\widetilde \Phi(s))ds
\end{equation*}
admits a unique solution $\widetilde \Phi\in C^0([0,T],B^1)$ and we have
\be
\label{diffduha}
\sup_{t\in[0,T]}\|\Phi(t)-\widetilde \Phi(t)\|_{B^1}\leq \eta.
\ee
\end{remark}

\subsection{Proof of Proposition \ref{propaver}}

This subsection is devoted to the proof of Proposition \ref{propaver}, which relies on a reformulation of the Duhamel formula for \fref{interm2}. 

\bs
\ni
{\em Step 1: reformulation of the Duhamel formula.} Introduce the following family of unbounded self-adjoint operators on $L^2(\RR^3)$
\be
\label{a}\forall \tau \in \RR,\qquad a(\tau)=-2B e^{i\tau H_{z}}ze^{-i\tau H_{z}}i\pa_x
\ee
with domain $B^2$. Note that, from \fref{H1} and Lemma \ref{lemBACM}, we deduce that, for all $\ell\in\NN$,
\be
\label{ftau}
\forall u\in B^2,\quad \forall \tau \in \RR,\qquad \|a(\tau)u\|_{L^2}\leq C\|u\|_{B^2}.
\ee
The Duhamel representation of \fref{interm2} reads
\be
\label{duha}
\widetilde \Phe(t)=\widetilde{\Psi_0}-\frac{i}{\eps}\int_0^t a\left(\frac{s}{\eps^2}\right)\widetilde \Phe(s)ds-i\int_0^t\left(-\Delta_{x,y}\widetilde\Phe(s)+ F_0\left(\widetilde\Phi^\eps(s)\right)\right)ds.
\ee
Introduce the primitive of $a$:
\be
\label{A}
\forall u\in B^2,\quad \forall \tau\in \RR,\qquad A(\tau)u=\int_0^\tau a(s)u\,ds,
\ee
which is well-defined as a Riemann integral, thanks to \fref{ftau}, and is such that
\be
\label{ftau2}
\forall u\in B^2,\quad \forall \tau \in \RR,\qquad \|A(\tau)u\|_{L^2}\leq C\tau\|u\|_{B^2}.
\ee
Now, we notice that if $\widetilde\Phe\in C^0([0,T],B^4)$, then by \fref{interm2} we have that  $\pa_t\widetilde\Phe\in C^0([0,T],B^2)$. Hence one can integrate by parts in the first integral of \fref{duha} and, if $m\geq 4$,  the following expression holds true for all $t\in [0,T]$, in the sense of functions in $C^0([0,T],L^2)$:
\bee
&\ds -\frac{i}{\eps}\int_0^t a\left(\frac{s}{\eps^2}\right)\widetilde \Phe(s)ds
&=i\eps\int_0^t A\left(\frac{s}{\eps^2}\right)\pa_t \widetilde \Phe(s)ds-i\eps  A\left(\frac{t}{\eps^2}\right)\widetilde \Phe(t)\\
&& =\int_0^t A\left(\frac{s}{\eps^2}\right)a\left(\frac{s}{\eps^2}\right)\widetilde\Phe(s)ds-i\eps  A\left(\frac{t}{\eps^2}\right)\widetilde \Phe(t)\\
&&\quad +\eps\int_0^t A\left(\frac{s}{\eps^2}\right)\left(-\Delta_{x,y}\widetilde\Phe(s)+ F_0\left(\widetilde\Phi^\eps(s)\right)\right)ds,
\eee
where we used \fref{interm2} to evaluate $i\pa_t \widetilde\Phe$. Finally, the Duhamel formula \fref{duha} becomes
\bea
&\widetilde \Phe(t)&=\widetilde{\Psi_0}+\int_0^t A\left(\frac{s}{\eps^2}\right)a\left(\frac{s}{\eps^2}\right)\widetilde \Phe(s)ds-i\eps  A\left(\frac{t}{\eps^2}\right)\widetilde \Phe(t)\nonumber\\
&& \quad +\eps\int_0^t A\left(\frac{s}{\eps^2}\right)\left(-\Delta_{x,y}\widetilde\Phe(s)+ F_0\left(\widetilde\Phi^\eps(s)\right)\right)ds\nonumber\\
&&\quad -i\int_0^t\left(-\Delta_{x,y}\widetilde\Phe(s)+ F_0\left(\widetilde\Phi^\eps(s)\right)\right)ds.\label{lesI}
\eea

\bs
\ni
{\em Step 2: approximation of the Duhamel formula.} Denote
$$\widehat \Phe(t)= \widetilde \Phe(t)+i\eps  A\left(\frac{t}{\eps^2}\right)\widetilde \Phe(t)$$
and rewrite \fref{lesI} as follows:
\bea
&\widehat \Phe(t)&=\widetilde{\Psi_0}+\int_0^t \left(A\left(\frac{s}{\eps^2}\right)a\left(\frac{s}{\eps^2}\right)+i\pa_x^2\right)\widetilde \Phe(s)ds-i\int_0^t\left(-\pa_y^2\widetilde\Phe(s)+ F_0\left(\widetilde\Phi^\eps(s)\right)\right)ds. \nonumber\\
&& \quad +\eps\int_0^t A\left(\frac{s}{\eps^2}\right)\left(-\Delta_{x,y}\widetilde\Phe(s)+ F_0\left(\widetilde\Phi^\eps(s)\right)\right)ds.\label{lesI2}
\eea
In this step, we prove that
\be
\label{e6}
\sup_{t\in[0,T]}\|\widetilde \Phe(t)-\widehat \Phe(t)\|_{B^1}\leq \eps\,C_M
\ee
and that
\be
\widehat \Phe(t)=\widetilde{\Psi_0}-i\int_0^t\left(A_0 \widehat \Phe(s)-\pa^2_y \widehat \Phe(s)+ F_0\left(\widehat\Phe(s)\right)+\eps f^\eps(s)\right)ds,\label{id2}
\ee
with 
\be
\sup_{t\in[0,T]}\left\|f^\eps\right\|_{B^1}\leq C_M.
\label{id3}
\ee
In order to prove this claim, we state two technical lemmas which are proved in the Appendix \ref{appC} so that the proof would be more readable.
\begin{lemma}
\label{lemtech2}
Let $V_c$ satisfy Assumptions \pref{confinement} and \pref{ass2}. Then, for all integer $\ell$, the operator $A(\tau)$ defined by \fref{A} satisfies
\be
\label{estiA}
\forall u\in C^0([0,T],B^{2n_{0}+\ell+8}),\qquad \sup_{t\in[0,T]}\left\|A\left(\frac{t}{\eps^2}\right)u(t)\right\|_{B^{\ell}}\leq C\|u\|_{C^0([0,T],B^{2n_{0}+\ell+8})},
\ee
where $n_0$ is as in Assumption \pref{ass2} and $C$ is independent of $\eps$.
\end{lemma}

\begin{lemma}
\label{lemtech1}
Let $V_c$ satisfy Assumptions \pref{confinement} and \pref{ass2}. Let $T>0$ and $m=4n_0+17$. Let $u\in C^0([0,T],B^m)$ such that $\pa_t u\in C^0([0,T],B^{m-2})$. Then we have, for all $\eps\in (0,1]$,
\be
\label{estimI1}
\sup_{t\in[0,T]}\left\|\int_0^t \left(A\left(\frac{s}{\eps^2}\right)a\left(\frac{s}{\eps^2}\right)+i\pa_x^2\right)u(s)ds+i\int_0^t A_0 u(s)ds\right\|_{B^1}\leq C\eps^2\|u\|
\ee
where $A_0$, $a$ and $A$ are respectively defined by \fref{A0}, \fref{a} and \fref{A} and where  $\|u\|$ denotes shortly $\|u\|_{C^0([0,T],B^m)}+\|\pa_tu\|_{C^0([0,T],B^{m-2})}$.
\end{lemma}

In order to apply these lemmas, we need some bounds for $\widetilde \Phi^\eps$ and $\pa_t\widetilde \Phi^\eps$. Let us fix $m=4n_0+17$, where $n_0$ is as in Assumption \ref{ass2} and assume that we have the uniform estimate \fref{estiaver1}. By \fref{derivation}, we deduce that
\be
\label{e2}
\|\Delta_{x,y}\widetilde \Phi^\eps\|_{C^0([0,T],B^{m-2})}
+\left\|e^{itH_{z}/\varepsilon^2}ze^{-itH_{z}/\varepsilon^2}\pa_x\widetilde\Phe\right\|_{C^0([0,T],B^{m-2})}\leq C_M.
\ee
Moreover, from \fref{tame}, we deduce that
\be
\label{e3}
\|F_0(\widetilde \Phi^\eps)\|_{C^0([0,T],B^m)}\leq C_M.
\ee
Hence, from \fref{interm2}, \fref{e2} and \fref{e3}, we get
\be
\label{e4}
\|\pa_t \widetilde \Phi^\eps\|_{C^0([0,T],B^{m-2})}\leq \frac{C_M}{\eps}.
\ee
Therefore, applying Lemmas \ref{lemtech2} and \ref{lemtech1} and using \fref{estiaver1}, \fref{e2}, \fref{e3} and \fref{e4} yield 
\be
\label{e5}
\sup_{t\in[0,T]}\left\|A\left(\frac{t}{\eps^2}\right)\widetilde \Phe(t)\right\|_{B^{2n_0+5}}\leq C_M,
\ee
\be
\label{e7}
\eps \sup_{t\in[0,T]}\left\|A\left(\frac{t}{\eps^2}\right)\left(-\Delta_{x,y}\widetilde\Phe(t)+ F_0\left(\widetilde\Phi^\eps(t)\right)\right)\right\|_{B^1}\leq \eps\,C_M
\ee
and
\be
\label{e8}
\sup_{t\in[0,T]}\left\|\int_0^t \left(A\left(\frac{s}{\eps^2}\right)a\left(\frac{s}{\eps^2}\right)+i\pa_x^2\right)\widetilde \Phe(s)ds+i\int_0^t A_0 \widetilde \Phe(s)ds\right\|_{B^1}\leq \eps\,C_M,
\ee
where we used that $m\geq 4n_0+17$, thus in particular $m\geq 4n_0+13$ and $m\geq 2n_0+11$.
Hence, from  \fref{e5}, we deduce \fref{e6} and
\be
\label{e9}
\|\pa_y^2(\widetilde\Phe-\widehat \Phe)\|_{C^0([0,T],B^1)}+\left\|F_0\left(\widetilde\Phe\right)-F_0\left(\widehat\Phe\right)\right\|_{C^0([0,T],B^1)}\leq \eps\,C_M,
\ee
where we also used the estimate \fref{estB1}. Moreover, from \fref{embed} and \fref{e5}, we get
\be
\label{e10}
\|A_0 (\widetilde\Phe-\widehat \Phe)\|_{C^0([0,T],B^1)}\leq \eps\,C_M.
\ee
Finally, inserting \fref{e7}, \fref{e8}, \fref{e9}, \fref{e10} in \fref{lesI2} yields \fref{id2} with the estimate \fref{id3}.

\bs
\ni
{\em Step 3: a stability result for the limit system.} First notice that \fref{id2} implies that $\widehat \Phe$ satisfies in the strong sense the equation
$$
i\pa_t \widehat \Phe=A_0\widehat \Phe-\pa^2_y\widehat \Phe+F_0(\widehat \Phe)+\eps f^\eps,\qquad \Phi(t=0)=\widetilde{\Psi_0}.
$$
which has the following mild formulation:
\be
\label{eq1}
\widehat \Phe(t)=e^{-it(A_0-\pa^2_y)}\widetilde {\Psi_0}-i\int_0^te^{-i(t-s)(A_0-\pa^2_y)}\left(F_0(\widehat \Phe(s))+\eps f^\eps\right)ds.
\ee
Apply now Proposition \ref{lemA0} {\em (ii)} with $\widetilde {\Psi_0}$ as initial data: there exists a maximal solution $\widetilde \Phi\in C^0([0,T_{max}),B^1)$ to the equation \fref{lim}. Assume that $T$ is such that $0<T<T_{max}$. Substracting \fref{lim} to \fref{eq1} leads, for all $t\leq T$, to
\bee
&\|\widehat \Phe(t)-\widetilde \Phi(t)\|_{B^1}&\leq \int_0^t\left\|F_0(\widehat \Phe(s))-F_0(\widetilde \Phi(s))\right\|_{B^1}ds+\eps \|f^\eps\|_{C^0([0,T],B^1)}\\
&&\leq C_M\left(\eps+\int_0^t\left\|\widehat \Phe(s)-\widetilde \Phi(s)\right\|_{B^1}ds\right),\\
\eee
where we used \fref{estB1}, \fref{id3} and $\|\widehat \Phe\|_{C^0([0,T],B^1)}\leq C_M$. Therefore, the Gronwall lemma gives the estimate \fref{estili} and the proof of Proposition \ref{propaver} is complete.
\qed

\section{Proof of the main theorem}
\label{lastsect}

This section is devoted to the proof of the main Theorem \ref{theo}.  Remark that the statement {\em (i)} is already proved in Proposition \ref{lemA0}. Let us prove the statement {\em (ii)} of Theorem \ref{theo}.

Let $\Psi_0\in B^1$. Denote by $\Psi^\eps\in C^0([0,+\infty),B^1)$ the solution of \fref{schrod}, \fref{cauchy}, \fref{poisson} and let $T_0\in (0,+\infty]$ be the maximal time given by Proposition \ref{prop1} {\em (i)}. We also introduce  the maximal solution $\Phi\in C^0([0,T_{max}),B^1)$ of the limit system \fref{eqA0}, given by Proposition \ref{lemA0}. Pick $T$ such that
$$0<T<\min(T_0,T_{max})$$
and let $\eta>0$.

Since $T<T_{max}$, according to Remark \ref{remstab}, one can define $\delta_{\eta/3,T}>0$ such that the following holds true. For all $\widetilde \Psi_0$ satisfying 
$$\|\Psi_0-\widetilde \Psi_0\|_{B^1}\leq \delta_{\eta/3,T},$$
 the equation \fref{lim} admits a unique solution $\widetilde \Psi\in C^0([0,T],B^1)$ and we have \fref{diffduha}:
\begin{equation*}
\sup_{t\in[0,T]}\|\Phi(t)-\widetilde \Phi(t)\|_{B^1}\leq \eta/3.
\end{equation*}
Next, we fix $m\geq 2$ according to Proposition \ref{propaver} and $\delta>0$ by
\be
\label{deltaaa}
\delta=\min\left(\frac{\eta}{3},\delta_{\eta/3,T}\right).
\ee
Since $T<T_0$, Proposition \ref{prop1} {\em (ii)} enables to choose $\widetilde{ \Psi_0}\in B^m$ and $\eps_\delta$ such that the corresponding solution $\widetilde \Psi^\eps$ of the intermediate system \fref{schrodI}, \fref{cauchyI}, \fref{poissonI} satisfies \fref{approx1} and \fref{approx2} for all $\eps\leq \eps_\delta$:
\be
\label{es2}
\|\Psi^\eps-\widetilde \Psi^\eps\|_{C^0([0,T],B^1)}\leq \delta\leq \frac{\eta}{3}
\ee
and $\widetilde \Psi^\eps$ is bounded in $C^0([0,T],B^m)$ uniformly with respect to $\eps$. 

Now, we remark that by \fref{es2} this initial data $\widetilde \Psi_0$ satisfies $$\|\Psi_0-\widetilde \Psi_0\|_{B^1}\leq \delta\leq \delta_{\eta/3,T}.$$ Hence, Remark \ref{remstab} gives that the solution $\widetilde \Phi$ of the equation \fref{lim} satisfies
$$
\|\Phi-\widetilde \Phi\|_{C^0([0,T],B^1)}\leq \frac{\eta}{3},
$$
or, equivalently,
\be
\label{es1}
\|e^{-itH_z/\eps^2}\Phi-e^{-itH_z/\eps^2}\widetilde \Phi\|_{C^0([0,T],B^1)}\leq \frac{\eta}{3},
\ee
Moreover, the uniform bound of $\widetilde \Psi^\eps$ in $C^0([0,T],B^m)$ enables to apply Proposition \ref{propaver}, which gives that the function $\widetilde \Pe$ satisfies
\be
\label{es3}
\|\widetilde \Pe-e^{-itH_z/\eps^2}\widetilde\Phi\|_{C^0([0,T],B^1)}\leq \delta\leq \frac{\eta}{3},
\ee
for $\eps\leq \eps_\delta$, where $\widetilde\Phi$ solves \fref{lim}.
Finally, \fref{es2}, \fref{es1} and \fref{es3} yield the existence of $\eps_0$ such that, for all $\eps\in (0,\eps_0]$ we have
\be
\label{es4}
\|\Pe-e^{-itH_z/\eps^2}\Phi\|_{C^0([0,T],B^1)}\leq \eta.
\ee
To conclude, it remains to remark that $T_0\geq T_{max}$. Indeed, if $T_0< T_{max}$, then we have, by Proposition \ref{prop1} {\em (i)},
$$\underset{\footnotesize \eps\to 0}{\mbox{\rm lim sup}}\|\Psi^\eps\|_{C^0([0,T_0],B^1)}=+\infty,$$
which implies by \fref{es4} that
$$\lim_{T\to T_0}\|\Phi(T)\|_{B^1}=+\infty.$$
This contradicts $T_0<T_{max}$. The proof of Theorem \ref{theo} is complete.
\qed

\begin{appendix}

\section{Proof of Lemma \ref{sobo2}}
\label{appZ}
First, by integrating by parts and applying Cauchy-Schwarz, we obtain
$$\|Bz\pa_xu\|_{L^2}^2=\int_{\RR^3} B^2z^2\,|\pa_xu|^2dxdydz=\int_{\RR^3} (B^2z^2u)(-\pa_x^2u)dxdydz\leq \|u\|_{B^2}^2.$$
Hence, the first properties stated in the Lemma are obvious from the definition \fref{HH}, and we shall only detail the proof of the equivalence of norms.

\bs
\ni
{\em Step 1: the case $m=1$}. From the definition \fref{HH} and the assumption \fref{H1} on $V_c$\,, we deduce that
\bee
&&\|H_\eps^{1/2}u\|_{L^2}^2=((-\pa^2_z+V_c)u,u)_{L^2}+\|(\eps\pa_x+iBz)u\|_{L^2}^2\\
&&=((-\pa^2_z+V_c)u,u)_{L^2}+B^2\|zu\|_{L^2}^2+\eps^2\|\pa_x u\|_{L^2}^2-2\eps B\im(zu,\pa_xu)_{L^2}\\
&&\geq \frac{1}{2}((-\pa^2_z+V_c)u,u)_{L^2}+ (\frac{a^2}{2}+B^2)\|zu\|_{L^2}^2+\eps^2\|\pa_x u\|_{L^2}^2-2\eps B\|zu\|_{L^2}\|\pa_x u\|_{L^2}^2\\
&&\geq \frac{1}{2}((-\pa^2_z+V_c)u,u)_{L^2}+ \frac{a^2}{4}\|zu\|_{L^2}^2+\frac{a^2}{a^2+4B^2}\eps^2\|\pa_x u\|_{L^2}^2\\
&&\geq C\|H_z^{1/2}u\|_{L^2}^2+C\eps^2\|\pa_xu\|_{L^2}^2.
\eee
Conversely, from \fref{equivB11} and \fref{HH}, we estimate directly
$$(H_\eps u,u)_{L^2}\leq C' \|H_z^{1/2}u\|_{L^2}^2+C'\eps^2\|\pa_xu\|_{L^2}^2.$$
For all $\eps\in (0,1]$, this yields the equivalence of norms \fref{Bell31}.

For $m\geq 2$, we will proceed by induction. For the clarity of the proof, let us introduce two notations. For $m\in \NN$, we denote by $(\mathcal P_m)$ the property\\[3mm]
\begin{tabular}{ll}$(\mathcal P_m)$: $\quad $ &there exists $\eps_m>0$ such that, for all $\eps\in (0,\eps_m]$ and for all $u\in B^m$, we have\\[2mm]
&$\ds \frac{1}{2}\|u\|_{B^m}^2\leq \|u\|_{L^2(\RR^3)}^2+\|\Delta_{x,y}^{m/2}u\|^2_{L^2(\RR^3)}+\|H_\eps^{m/2}u\|^2_{L^2(\RR^3)}\leq 2\|u\|_{B^m}^2,$\end{tabular}\\[3mm]
and by $(\mathcal Q_m)$ the property\\[3mm]
\begin{tabular}{ll}$(\mathcal Q_m)$:$\quad$& there exists $C_m>0$ such that, for all $u\in B^m$ and $\eps\in (0,1]$,\\[2mm]& the operator $A_m=\frac{1}{\eps}(H_\eps^m-H_z^m)$ satisfies  $|(A_mu,u)_{L^2}|\leq C_m\|u\|_{B^m}^2$.\end{tabular}\\[3mm]
Note that the lemma will proved if we show that $(\mathcal P_m)$ holds true for all $m\geq 0$. Note also that, up to a possible modification of the sequence $(\eps_m)_{m\in \NN}$, this sequence can be chosen nonincreasing.

\bs
\ni
{\em Step 2: $(\mathcal Q_m)$ implies $(\mathcal P_m)$}. Let $m\geq 0$ be fixed. From $(\mathcal Q_m)$, we deduce that
\bee
&\|H_\eps^{m/2}u\|_{L^2}^2=(H_\eps^mu,u)_{L^2}&=(H_z^mu,u)_{L^2}+\eps (A_mu,u)_{L^2}\\
&&=\|H_z^{m/2}u\|_{L^2}^2+\eps (A_mu,u)_{L^2}\,,
\eee
thus
\be
\label{normin}
\|H_z^{m/2}u\|_{L^2}^2-\eps C_m\|u\|^2_{B^m}\leq \|H_\eps^{m/2}u\|^2_{L^2}\leq \|H_z^{m/2}u\|_{L^2}^2+\eps C_m\|u\|^2_{B^m}\,.
\ee
Setting
$$\eps_m=\frac{1}{2C_m}\,,$$
 we deduce directly from \fref{Bell} and \fref{normin} that, for $\eps\leq\eps_m$, 
$$\frac{1}{2}\|u\|_{B^m}^2\leq \|u\|_{L^2(\RR^3)}^2+\|\Delta_{x,y}^{m/2}u\|^2_{L^2(\RR^3)}+\|H_\eps^{m/2}u\|^2_{L^2(\RR^3)}\leq 2\|u\|_{B^m}^2.$$
We have proved $(\mathcal P_m)$.

\bs
\ni
{\em Step 3: proof of $(\mathcal Q_m)$ for $m=0$ and $1$.} For $m=0$, choose $A_0=0$ and  $(\mathcal Q_0)$ is obvious. Let us prove  $(\mathcal Q_1)$. From \fref{HH}, we have
\be
\label{A1}
H_\eps=H_z+\eps A_1,\quad \mbox{with}\quad A_1=-2iBz\pa_x-\eps\pa^2_x.
\ee
For all $u\in B^1$, we have
$$|(A_1u,u)_{L^2}|=|-2iB(\pa_xu,zu)_{L^2}+\eps \|\pa_x u\|_{L^2}^2|\leq C(\|zu\|_{L^2}^2+\|\pa_xu\|_{L^2}^2)\leq C_1\|u\|_{B^1}^2,$$ 
where we applied Cauchy-Schwarz and Lemma \ref{lemBACM}. We have proved $(\mathcal Q_1)$.

\bs
\ni
{\em Step 4: proof of $(\mathcal Q_m)$ for $m\geq 2$.} We shall now proceed by  induction. Let $m\geq 2$ and assume that $(\mathcal Q_{m-2})$ and $(\mathcal Q_{m-1})$ hold true. Let us prove  $(\mathcal Q_m)$. We compute
\bee
&H_\eps^{m}&=(H_z+\eps A_1)H_\eps^{m-2}(H_z+\eps A_1)\\
&&=H_zH_\eps^{m-2}H_z+\eps A_1 H_\eps^{m-1}+\eps H_\eps^{m-1}A_1\\
&&=H_z^m+\eps H_zA_{m-2}H_z+\eps A_1 H_\eps^{m-1}+\eps H_\eps^{m-1}A_1
\eee
where we have applied $(\mathcal Q_{m-2})$. Hence, denoting
\be
\label{Am}
A_m=H_zA_{m-2}H_z+ A_1 H_\eps^{m-1}+ H_\eps^{m-1}A_1\,,
\ee
we obtain
$$H_\eps^m=H_z^m+\eps A_m$$
and, for all $u\in B^m$, we get from the definition \fref{Am} that
$$
|(A_mu,u)_{L^2}|\leq |(H_zA_{m-2}H_z u,u)_{L^2}|+2|(H_\eps^{m-1}u,A_1u)|_{L^2}\,,
$$
where we used that $H_\eps^{m-1}$ and the operator $A_1$ defined by \fref{A1} are selfadjoint. It remains to estimate the two terms in the right-hand side of this inequality. The first one can be estimated as follows:
$$
|(H_zA_{m-2}H_z u,u)_{L^2}|=|(A_{m-2}H_zu, H_zu)_{L^2}|\leq C_{m-2}\|H_zu\|_{B^{m-2}}^2 \leq C_{m-2}\|u\|_{B^m}^2\,,
$$
where we used $(\mathcal Q_{m-2})$ and \fref{derivation}. The second one can be estimated as follows:
\bee
&|(H_\eps^{m-1}u,A_1u)|_{L^2}&=\left|\left(H_\eps^{m-1}u,(i\pa_x)(-2Bzu+i\pa_xu)\right)\right|_{L^2}\\
&&=\left|\left(H_\eps^{\frac{m-1}{2}}(i\pa_x u),H_\eps^{\frac{m-1}{2}}(-2Bzu+i\pa_xu)\right)\right|_{L^2}\\&&\leq \left\|H_\eps^{\frac{m-1}{2}}(i\pa_x u)\right\|_{L^2}\left\|H_\eps^{\frac{m-1}{2}}(-2Bzu+i\pa_xu)\right\|_{L^2}\\
&&\leq C\|\pa_x u\|_{B^{m-1}}\|z u\|_{B^{m-1}}+C\|\pa_x u\|_{B^{m-1}}^2\\
&&\leq C\|u\|_{B^m}^2\,,
\eee
where we used that $H_\eps$ commutes with $\pa_x$, the  Cauchy-Schwarz inequality, the property $(\mathcal P_{m-1})$ and, at the last step, \fref{derivation}. Therefore, we have proved that
$$
|(A_mu,u)_{L^2}|\leq C_m \|u\|_{B^m}^2\,,
$$
which proves $(\mathcal Q_m)$. The proof of the lemma is complete.

\section{Proof of Lemma \ref{tame estimate}}
\label{appA}

For readability, we introduce in this appendix the following notation:
$$\forall (x,y,z)\in \mathbb{R}^3,\ \forall \alpha\in \{0,1\},\ \forall \eps\in (0,1),\ \ r^{\eps}_{\alpha}(x,y,z)=\sqrt{x^2+y^2+\alpha\eps^2z^2}.$$
With this notation, for all $u\in B^1$, and $\alpha\in \{0,1\}$, the nonlinearity $F_{\alpha}(u)$ defined in \fref{unif} reads
$$F_{\alpha}(u)=\left(\frac{1}{4\pi r_{\alpha}^{\eps}}\ast(|u|^2)\right)u.$$ 

In order to prove the estimates stated in Lemma \ref{tame estimate}, we prove the following technical lemma on the Poisson nonlinearity.

\begin{lemma}
\label{lempoisson} The following estimates hold.\\
(i) There exists a positive constant $C$ that does not depend on $\eps\in (0,1]$ or $\alpha\in \{0,1\}$ such that 
\be
\label{i}
 \forall u,v\in H^1(\mathbb{R}^3),\ \left\|\frac{1}{r_{\alpha}^{\eps}}\ast \left(uv\right)\right\|_{L^{\infty}}\leq C \|u\|_{H^1}\|v\|_{H^1}.
 \ee 
(ii) There exists a positive constant $C$ that does not depend on $\eps\in (0,1]$ or $\alpha\in \{0,1\}$ such that, if $D$ denotes a derivative with respect to $x,y$ or $z$, \be
\label{ii}
\forall u,v\in H^1(\RR^3),\ \forall v\in H^1(\RR^3),\ \left\|D\left(\frac{1}{r_{\alpha}^{\eps}}\ast \left(uv\right)\right)\right\|_{L^3_{x,y}L^{\infty}_{z}}\leq C \|u\|_{H^1}\|v\|_{H^1}
\ee
(iii) For any integer $k$, let $\beta=(\beta_x,\beta_y,\beta_z)\in \NN^3$ be a multiinteger of length $|\beta|=\beta_x+\beta_y+\beta_z=k$ and let  $D^\beta=\pa_x^{\beta_x}\pa_y^{\beta_y}\pa_z^{\beta_z}$ be the associated derivative. Then there exists a positive constant $C_{k}$ depending only on $k$ such that
\be
\label{iii}
\forall u\in H^k,\ \left\|D^\beta\left(\frac{1}{r_{\alpha}^{\eps}}\ast |u|^2\right)\right\|_{L^3_{x,y}L^{\infty}_{z}}\leq C_{k} \|u\|_{H^1}\|u\|_{H^k}.
\ee
\end{lemma}
\begin{proof}
Noting that, for all $(x,y)\in \RR^2$, 
\begin{equation}
\label{first}
\left\|\left(\frac{1}{r^{\eps}_{\alpha}}\ast \left(uv\right)\right)(x,y,\cdot)\right\|_{L^{\infty}(\RR)}\leq \int_{\RR^2}\frac{1}{\sqrt{(x-x')^2+(y-y')^2}}\left\|uv(x',y',\cdot)\right\|_{L^1(\RR)}dx'dy',
\end{equation}
we only need estimates for the convolution with $\frac{1}{\sqrt{x^2+y^2}}$ in $\RR^2$. Here, we refer the reader to Lemma B.1 of \cite{BAMP} where it was shown that for any $f\in L^p(\RR^2)\cap L^1(\RR^2)$ with $2<p\leq \infty$, the following bound holds:
\be
\label{second}
\left\|\frac{1}{\sqrt{x^2+y^2}}\ast f\right\|_{L^{\infty}(\RR^2)}\leq C_{p}\|f\|^{\theta}_{L^p(\RR^2)}\|f\|^{1-\theta}_{L^1(\RR^2)}
\ee 
where $\theta=p/(2p-2)$. Moreover, from Cauchy-Schwarz and Sobolev embeddings, we deduce that for all $p\in [1,+\infty)$,
\bee
&\left\|\left\|uv(x,y,\cdot)\right\|_{L^1(\RR)}\right\|_{L^p(\RR^2)}
&\leq \left\|\left\|u(x,y,\cdot)\right\|_{L^2(\RR)}\left\|v(x,y,\cdot)\right\|_{L^2(\RR)}\right\|_{L^p(\RR^2)}\\
&&\leq \|u\|_{L^{2p}_{x,y}L^2_z}\|v\|_{L^{2p}_{x,y}L^2_z}\leq \|u\|_{H^1(\RR^3)}\|v\|_{H^1(\RR^3)}.
\eee Combined with \fref{first} and \fref{second}, this proves Item {\it (i)}. 

In order to prove Item {\it (ii)}, consider a first order derivative $D$ with respect to $x,y$ or $z$ and let $u,v\in H^1(\RR^3)$. Usual properties of the convolution give
$$D\left(\frac{1}{r_{\alpha}^{\eps}}\ast \left(uv\right)\right)=\frac{1}{r_{\alpha}^{\eps}}\ast D\left(uv\right)=\frac{1}{r_{\alpha}^{\eps}}\ast \left(D(u)v+uD(v)\right).$$
Using \fref{first} combined with the generalized Young formula gives
\bea
\label{oups}
\nonumber \left\|\frac{1}{r_{\alpha}^{\eps}}\ast \left(D(u)v+uD(v)\right)\right\|_{L^3_{x,y}L^{\infty}_{z}}&\leq& \left\|\frac{1}{ \sqrt{x^2+y^2}}\ast \left\|D(u)v+uD(v)\right\|_{L_{z}^1}\right\|_{L^3_{x,y}}\\
&&\leq C \left\|\left\|D(u)v+uD(v)\right\|_{L_{z}^1}\right\|_{L^{6/5}_{x,y}}
\eea
since the function $x\mapsto \frac{1}{\sqrt{x^2+y^2}}$ belongs to $L_{w}^2(\RR^2)$. We end the proof of Item {\it (ii)} noting that, thanks to Sobolev embeddings,
\bee
\left\|\left\|D(u)v+uD(v)\right\|_{L_{z}^1}\right\|_{L^{6/5}_{x,y}}&\leq& C\|D(u)\|_{L^2}\|v\|_{L^3_{x,y}L^2_{z}}+C\|D(v)\|_{L^2}\|u\|_{L^3_{x,y}L^2_{z}}\\&\leq& C \|u\|_{H^1}\|v\|_{H^1}.
\eee

In order to prove Item {\it (iii)}, we follow the same lines with derivatives of higher orders. Consider the derivative $D^\beta$ where $\beta=(\beta_x,\beta_y,\beta_z)\in \NN^3$ is a multiinteger of length $|\beta|=\beta_x+\beta_y+\beta_z=k$. Usual properties of the convolution gives
\begin{equation*}
\left\|D^\beta\left(\frac{1}{ r_{\alpha}^{\eps}}\ast \left(|u|^2\right)\right)\right\|_{L^3_{x,y}L^{\infty}_{z}}=\left\|\frac{1}{r_{\alpha}^{\eps}}\ast D^\beta\left(|u|^2\right)\right\|_{L^3_{x,y}L^{\infty}_{z}}.
\end{equation*}
Again, using \fref{first} combined to the generalized Young's formula lead to:
\be
\label{four2}
\left\|\frac{1}{r_{\alpha}^{\eps}}\ast D^\beta\left(|u|^2\right)\right\|_{L^3_{x,y}L^{\infty}_{z}}\leq \left\|\frac{1}{\sqrt{x^2+y^2}}\ast \left\|D^\beta\left(|u|^2\right)\right\|_{L_{z}^1(\RR)}\right\|_{L^3_{x,y}}\leq C\left\|D^\beta\left(|u|^2\right)\right\|_{L^{6/5}_{x,y}L^1_{z}}.
\ee
We now write
$$D^\beta(u\overline{u})=\sum_{\beta'\leq \beta} C_{\beta'}D^{\beta'}(u)D^{\beta-\beta'}(\overline{u}), $$
where the sum is over the set of multiintegers $\beta'=(\beta'_x,\beta'_y,\beta'_z)$ such that $\beta'_x\leq \beta_x$, $\beta'_y\leq \beta_y$ and $\beta'_z\leq \beta_z$. Thus, combining \fref{four2} with Sobolev embeddings gives as above
\bee
\left\|D^\beta\left(\frac{1}{r^{\eps}_{\alpha}}\ast \left(|u|^2\right)\right)\right\|_{L^3_{x,y}L^{\infty}_{z}}&\leq&
 C\sum_{|\beta'|=k}\|D^{\beta'}(u)\|_{L^2}\|u\|_{L^3_{x,y}L^2_z}\\
 &&+C\sum_{|\beta'|=\ell,\, 0\leq \ell<k}\|D^{\beta'}u\|_{L^3_{x,y}L^2_z}\|D^{\beta-\beta'}(\overline{u})\|_{L^2}\\
 &&\leq C\|u\|_{H^{k}}\|u\|_{H^{1}}+C\sum_{\ell=0}^{k-1}\|u\|_{H^{\ell+1}}\|u\|_{H^{k-\ell}}.
 \eee
We conclude noting that, by interpolation, for all $\ell\leq k-1$, 
\be
\label{interp}
\|u\|_{H^{\ell+1}}\|u\|_{H^{k-\ell}}\leq \|u\|_{H^{k}}\|u\|_{H^{1}}.
\ee
\end{proof}
\ni
{\em Proof of Lemma \ref{tame estimate}.} We first prove \fref{estB1}. In that view, let us fix $u$ and $v$ in $B^1$, $\eps\in(0,1)$ and $\alpha\in \{0,1\}$, and note that
\be
\label{rewrite}
\hskip -0.15cm F_{\alpha}(u)-F_{\alpha}(v)=\left(\frac{1}{4\pi r_{\alpha}^{\eps}}\ast\big[(|u|+|v|)(|u|-|v|)\big]\right)u+\left(\frac{1}{4\pi r_{\alpha}^{\eps}}\ast\left(|v|^2\right)\right)\left(u-v\right)
\ee
Hence \fref{estB1} is a straightforward consequence of the following claim. There exists a positive constant $C$ such that for all $u_{1},u_{2}$ and $u_{3}\in B^1$
\be
\label{u1u2u3}
\left\|\left(\frac{1}{r_{\alpha}^{\eps}}\ast(u_{1}u_{2})\right)u_{3}\right\|_{B^1}\leq C\|u_{1}\|_{B^1}\|u_{2}\|_{B^1}\|u_{3}\|_{B^1}.
\ee
{\em Proof of the claim \fref{u1u2u3}.} According to Lemma \ref{lemBACM} we have
\begin{multline}
\label{claim1}
\left\|\left(\frac{1}{r_{\alpha}^{\eps}}\ast(u_{1}u_{2})\right)u_{3}\right\|_{B^1}\leq C\left\|\left(\frac{1}{r_{\alpha}^{\eps}}\ast(u_{1}u_{2})\right)u_{3}\right\|_{H^1}+C\left\|\sqrt{V_{c}}\left(\frac{1}{r_{\alpha}^{\eps}}\ast(u_{1}u_{2})\right)u_{3}\right\|_{L^2}
\end{multline}
First, applying \fref{i} and then Lemma \ref{lemBACM},
\be
\label{claim2}
\left\|\left(\frac{1}{r_{\alpha}^{\eps}}\ast(u_{1}u_{2})\right)u_{3}\right\|_{L^2}\leq \left\|\frac{1}{r_{\alpha}^{\eps}}\ast(u_{1}u_{2})\right\|_{L^{\infty}}\|u_{3}\|_{L^2}\leq C\|u_{1}\|_{B^1}\|u_{2}\|_{B^1}\|u_{3}\|_{B^1}.
\ee
Similarly, we have
\bea
\label{claim3}
\nonumber\left\|\sqrt{V_{c}}\left(\frac{1}{r_{\alpha}^{\eps}}\ast(u_{1}u_{2})\right)u_{3}\right\|_{L^2}&&\leq \left\|\frac{1}{r_{\alpha}^{\eps}}\ast(u_{1}u_{2})\right\|_{L^{\infty}}\|\sqrt{V_{c}}u_{3}\|_{L^2}\\
&&\leq C\|u_{1}\|_{B^1}\|u_{2}\|_{B^1}\|u_{3}\|_{B^1}
\eea
Moreover, if $D$ denotes any differential operator of order $1$,
\be
\label{claimdec}
D\left(\left(\frac{1}{r_{\alpha}^{\eps}}\ast(u_{1}u_{2})\right)u_{3}\right)=D\left(\frac{1}{r_{\alpha}^{\eps}}\ast (u_{1}u_{2})\right)u_{3}+\left(\frac{1}{r_{\alpha}^{\eps}}\ast(u_{1}u_{2})\right)D(u_{3}).
\ee
Applying the H\"older inequality leads to
\bea
\left\|D\left(\frac{1}{r_{\alpha}^{\eps}}\ast (u_{1}u_{2})\right)u_{3}\right\|_{L^2}&\leq& \left\|D\left(\frac{1}{r_{\alpha}^{\eps}}\ast (u_{1}u_{2})\right)\right\|_{L^3_{x,y}L^{\infty}_{z}}\|u_{3}\|_{L^6_{x,y}L_{z}^2}\nonumber \\
&\leq& C\|u_{1}\|_{B^1}\|u_{2}\|_{B^1}\|u_{3}\|_{B^1}\,,\label{claim4}
\eea
where we used \fref{ii}. Finally, using \fref{i},
\bea
&\left\|\left(\ds \frac{1}{r_{\alpha}^{\eps}}\ast(u_{1}u_{2})\right)D(u_{3})\right\|_{L^2}
&\leq  \left\|\frac{1}{r_{\alpha}^{\eps}}\ast(u_{1}u_{2})\right\|_{L^\infty}\left\|D(u_{3})\right\|_{L^2}\nonumber\\
&&\leq C\|u_{1}\|_{B^1}\|u_{2}\|_{B^1}\|u_{3}\|_{B^1}\,.\label{claim5}
\eea
We deduce the claim \fref{u1u2u3} by combining \fref{claim1} with \fref{claim2}, \fref{claim3}, \fref{claimdec}, \fref{claim4} and \fref{claim5}. 

In order to prove \fref{tame}, consider a positive integer $m$ and fix $u\in B^m$. According to Lemma \ref{lemBACM}, we only need to estimate $\|F_{\alpha}(u)\|_{H^m}$ and $\|V_{c}^{m/2}F_{\alpha}(u)\|_{L^2}$. In that view, we readily have
\bea
&\left\|(1+V_{c}^{m/2})F_{\alpha}(u)\right\|_{L^2(\RR^3)}
&\leq\left\|\frac{1}{4\pi r_{\alpha}^{\eps}}\ast |u|^2\right\|_{L^{\infty}(\RR^3)}\|(1+V_{c}^{m/2})u\|_{L^2(\RR^3)}\nonumber\\
&&
\leq C \|u\|^2_{B^1(\RR^3)}\|u\|_{B^m}
\label{tame1}
\eea
where we applied \fref{i} and Lemma \ref{lemBACM}.

Now, let $D^\beta$ denote any derivative of length $m$ and write 
$$D^\beta(F_{\alpha}(u))=\sum_{\beta'\leq \beta}C_{\beta'}D^{\beta'}\left(\frac{1}{4\pi r_{\alpha}^{\eps}}\ast |u|^2\right)D^{\beta-\beta'}(u).$$
Hence,
\bea
\nonumber \left\|D^\beta F_{\alpha}(u)\right\|_{L^2(\RR^3)}&\leq& C\sum_{|\beta'|=m}\left\|\frac{1}{4\pi r_{\alpha}^{\eps}}\ast |u|^2\right\|_{L^{\infty}}\|D^{\beta'}u\|_{L^2}\\
&&\hspace*{-2cm}+C\sum_{1\leq |\beta'|\leq m}\left\|D^{\beta'}\left(\frac{1}{4\pi r_{\alpha}^{\eps}}\ast |u|^2\right)\right\|_{L^3_{x,y}L^{\infty}_{z}}\|D^{\beta-\beta'}(u)\|_{L^6_{x,y}L^2_{z}}\\
\nonumber &&\hspace*{-2cm}\leq C\|u\|_{H^1}^2\|u\|_{H^{m}}+C\sum_{\ell=1}^m\|u\|_{H^1}\|u\|_{H^\ell}\|u\|_{H^{m-\ell+1}}
\eea
where we applied \fref{i}, \fref{iii} and Sobolev embeddings. Using the interpolation estimate \fref{interp} gives
\be
\label{tame3}
 \left\|D^mF_{\alpha}(u)\right\|_{L^2(\RR^3)}\leq C\|u\|^2_{H^1}\|u\|_{H^m}\leq C\|u\|^2_{B^1}\|u\|_{B^m}\,.
\ee
We conclude the proof of \fref{tame} combining \fref{tame1} and \fref{tame3}. This ends the proof of Lemma \ref{tame estimate}.
\qed

\section{Proof of Lemma \ref{Poisson}}
\label{appB}
\ni
In this section, we set for simplicity $X=(x,y)\in \RR^2$. In order to prove estimate \fref{diffnoyau}, we first study the difference between both convolution kernels.\\
{\it First Step: Difference between the convolution kernels}\\
Let $u,v$ be two functions of $B^2$. Denote
\be
\label{deltadef}
\delta(u,v)(X,z)=\intd \ints \left(\frac{1}{\sqrt{|X-X'|^2+\eps^2(z-z')^2}}-\frac{1}{|X-X'|}\right)u(X',z')v(X',z')dz'dX'.
\ee
We split the integral as follows:
$$\delta(u,v)(X,z)=\delta^+(u,v)(X,z)+\delta^-(u,v)(X,z)=\int_{X'\in\Omega^+} \int_{z'\in\RR}+\int_{X'\in\Omega^-} \int_{z'\in\RR},$$
where
$$\Omega^+=\{X'\in \RR^2,\ |X-X'|>\eps\},\qquad \Omega^-=\{X'\in \RR^2,\ |X-X'|<\eps\}.$$
For all $\eta,\mu\in \RR$, and $X'\neq X$, we have
\be
\label{integrale}
\frac{1}{\sqrt{|X-X'|^2+\eps^2\eta^2}}-\frac{1}{\sqrt{|X-X'|^2+\eps^2\mu^2}}=\int_{\eps\mu}^{\eps\eta}\frac{-\xi}{\left(|X-X'|^2+\xi^2\right)^{3/2}}d\xi
\ee 
and
\be
\label{idiot}
\frac{1}{\sqrt{|X-X'|^2+\eps^2\eta^2}}-\frac{1}{\sqrt{|X-X'|^2+\eps^2\mu^2}}\leq \frac{2}{|X-X'|},
\ee
Besides, a simple study gives  
\be
\label{malin}
\forall X'\neq X,\ \forall \xi\in\RR,\ \frac{|\xi|}{\left(|X-X'|^2+\xi^2\right)^{3/2}}\leq \frac{2}{3\sqrt{3}}\frac{1}{|X-X'|^2}.
\ee
Equation \fref{integrale}, combined with \fref{idiot} and \fref{malin} allows us to claim that for all $\theta \in (0,1)$, 
\be
\label{estintegrale}
\left|\frac{1}{\sqrt{|X-X'|^2+\eps^2\eta^2}}-\frac{1}{\sqrt{|X-X'|^2+\eps^2\mu^2}}\right|\\
\leq C \eps^{\theta}|\eta-\mu|^{\theta}\frac{1}{|X-X'|^{1+\theta}}.
\ee
Now, applying \fref{estintegrale} with $\eta=z-z'$, $\mu=z'$ and $\theta=3/8$ leads to
\bee
\nonumber &&\left|\int_{\Omega^+}\int_{\RR}\left(\frac{1}{\sqrt{|X-X'|^2+\eps^2(z-z')^2}}-\frac{1}{\sqrt{|X-X'|^2+\eps^2z'^2}}\right)u(X',z')v(X',z')dz'dX'\right|\\
\nonumber && \hskip 1cm \leq C \eps^{3/8} |z|^{3/8}\int_{\Omega^+}\frac{1}{|X-X'|^{11/8}}\|u(X',\cdot)v(X',\cdot)\|_{L^1_{\RR}}dX'\\
&& \hskip 1 cm \leq C \eps^{3/8}|z|^{3/8}\frac{1}{\eps^{1/24}}\|u\|_{L^6_{X}L^2_{z}}\|v\|_{L^6_{X}L^2_{z}}
\leq C \eps^{1/3}|z|^{3/8}\|u\|_{B^1}\|v\|_{B^1}.
\eee
where we used the H\"older inequality and Sobolev embeddings. Similarly, applying \fref{estintegrale} with $\eta=z'$ , $\mu=0$ and $\theta=3/4$ leads to
\bee
\nonumber \left|\int_{\Omega^+}\int_{\RR}\left(\frac{1}{\sqrt{|X-X'|^2+\eps^2z'^2}}-\frac{1}{|X-X'|}\right)|u(X',z')||v(X',z')|dz'dX'\right|\\ 
\nonumber && \hskip -8cm \leq C \eps^{3/4} \int_{\Omega^+}\frac{1}{|X-X'|^{7/4}}\||z'|^{3/4}u(X',\cdot)v(X',\cdot)\|_{L^1_{\RR}}dX'\\
\nonumber && \hskip -8cm \leq C \eps^{3/4}\frac{1}{\eps^{5/12}}\|z^{3/8}u\|_{L^6_{X}L^2_{z}}\|v\|_{L^6_{X}L^2_{z}}\\&& \hskip -8cm \leq C \eps^{1/3}\|u\|_{B^2}\|v\|_{B^1}.
\eee
We have proved that
\be
\label{inequal1}
|\delta^+(u,v)(X,z)|\leq C \eps^{1/3}(1+|z|^{3/8})\|u\|_{B^2}\|v\|_{B^1} 
\ee
Consider now $\delta^-$. Using \fref{integrale} again leads to
\be
\label{remalin}
|\delta^-(u,v)(X,z)|\leq\int_{\Omega^-}\int_{\RR}\int_{\RR}\frac{|\xi|}{(|X-X'|^2+\xi^2)^{3/2}}|u(X',z')||v(X',z')|d\xi dz'dX'.
\ee
Moreover, a simple computation gives
$$\int_{\RR}\frac{|\xi|}{(|X-X'|^2+\xi^2)^{3/2}}d\xi=\frac{2}{|X-X'|}.$$
Hence, \fref{remalin} gives
\bea
\label{d-}
\nonumber |\delta^-(u,v)(X,z)|&\leq& C \int_{\Omega-}\int_{\RR}\frac{1}{|X-X'|}|u(X',z')||v(X',z')|dz'dX'\\ 
\nonumber &\leq& C\eps^{1/3}\|u\|_{L^6_{X}L^2_{z}}\|v\|_{L^6_{X}L^2_{z}}\\
&\leq& C \eps^{1/3}\|u\|_{B^1}\|v\|_{B^1}.
\eea
Combining \fref{inequal1} and \fref{d-} allows to conclude that 
\be
\label{final}
|\delta(u,v)(X,z)|\leq C\eps^{1/3}(1+\sqrt{V_c(z)})\|u\|_{B^2}\|v\|_{B^1},
\ee
where we have used $z^{3/8}\leq C(1+\sqrt{V_c(z)})$, deduced from \fref{H1}.

\bs
\ni
{\it Step 2: Difference between the nonlinearities}. In order to prove Lemma \ref{Poisson}, we need to estimate the following quantity in $B^1$:
\be
\label{dec}
F_{1}(u)-F_{0}(u)=\delta(u,\overline{u})\,u,\ee
where $u\in B^2$ is given. According to Lemma \ref{lemBACM}, we have
$$\|F_{1}(u)-F_{0}(u)\|_{B^1}\leq C\|\sqrt{V_{c}}\left(F_{1}(u)-F_{0}(u)\right)\|_{L^2} +C\|F_{1}(u)-F_{0}(u)\|_{H^1}.$$
First, we deduce from \fref{final} that
\be
\label{P1}
\|(1+\sqrt{V_c})\delta(u,\overline{u})\,u\|_{L^2}\leq C\eps^{1/3}\|(1+V_c)u\|_{L^2}\|u\|_{B^2}\|v\|_{B^1}\leq C\eps^{1/3}\|u\|_{B^2}^3,
\ee
where we used Lemma \ref{lemBACM}. Let now $D$ denote a first order derivative with respect to $x$, $y$ or $z$. We clearly have
\bea
\label{dedec}
\nonumber &&\|D\left(F_{1}(u)-F_{0}(u)\right)\|_{L^2}\\
&&\qquad \leq\left\|\left(\frac{1}{\sqrt{|X|^2+\eps^2z^2}}-\frac{1}{|X|}\right)\ast\left(D(u)\overline{u}+uD(\overline{u})\right)u\right\|_{L^2}+\|\delta(u,\overline{u})D(u)\|_{L^2}\nonumber\\
&&\qquad \leq 2|\delta(\overline{u},D(u))u\|_{L^2}+\|\delta(u,\overline{u})D(u)\|_{L^2}.
\eea
According to \fref{final}, we have
\be
\label{ggauche}
\|\delta(\overline{u},D(u))u\|_{L^2}\leq C\eps^{1/3}\|(1+\sqrt{V_c})u\|_{L^2}\|u\|_{B^2}\|D(u)\|_{B^1}
\leq C\eps^{1/3}\|u\|^3_{B^2}
\ee
and
\be
\label{ddroite}
\|\delta(u,\overline{u}) D(u)\|_{L^2}\leq C\eps^{1/3}\|(1+\sqrt{V_c})D(u)\|_{L^2}\|u\|_{B^2}\|u\|_{B^1}\leq C\eps^{1/3}\|u\|^3_{B^2},
\ee
where we used again Lemma \ref{lemBACM}. Combining \fref{dec}, \fref{P1}, \fref{dedec},\fref{ggauche} and \fref{ddroite} gives \fref{diffnoyau}. The proof of Lemma \ref{Poisson} is complete.
\qed

\section{Proof of the technical Lemmas \ref{lemtech2} and \ref{lemtech1}}
\label{appC}

Let us develop the operators $a$ and $A$ defined by \fref{a} and \fref{A} on the eigenbasis $\chi_p$. We have
$$
a(\tau)u=-\sum_{p\geq 0}\sum_{q\geq 0}e^{i\tau (E_p-E_q)}a_{pq}\,i\pa_x u_q\,\chi_p
$$
where we have introduced the coefficients
\be
\label{apq}
a_{pq}=\left\langle 2Bz \chi_p \chi_q\right\rangle.
\ee
Recall that, by Assumption \ref{confinement}, the potential $V_c$ is even, so for all $p$, the function $(\chi_p(z))^2$ is even. Therefore, we have
$$
\forall p\in \NN,\quad a_{pp}=\left\langle 2Bz \chi_p^2 \right\rangle=0,
$$
thus
\be
\label{a2}
a(\tau)u=-\sum_{p\geq 0}\sum_{q\neq p}e^{i\tau (E_p-E_q)}a_{pq}\,i\pa_x u_q\,\chi_p\,.
\ee
Let us now integrate this formula in order to compute the operator $A$ defined by \fref{A}:
\be
\label{A2}
A(\tau)u=i\sum_{p\geq 0}\sum_{q\neq p}\frac{e^{i\tau (E_p-E_q)}-1}{E_p-E_q}a_{pq}\,i\pa_x u_q\,\chi_p\,.
\ee
Before proving Lemmas \ref{lemtech2} and \ref{lemtech1}, let us give a useful estimate on coefficients $a_{pq}$. For all $p\in \NN$, $q\in \NN$, $k\in \NN$ we have
\be
\label{borneapq}
|a_{pq}|\leq C\frac{E_q^{(k+1)/2}}{E_p^{k/2} }.
\ee
Indeed, we have
\bee
&\left|E_p^{k/2}a_{pq}\right|=2B\left|\left(H_z^{k/2}\chi_p,z\chi_q\right)_{L^2}\right|&=2B\left|\left(\chi_p,H_z^{k/2}(z\chi_q)\right)_{L^2}\right|\\
&&\leq 2B\|H_z^{k/2}(z\chi_q)\|_{L^2}\\
&&\leq 2B\|z\chi_q\|_{B^k}\\
&&\leq C\|\chi_q\|_{B^{k+1}}\leq CE_q^{(k+1)/2},
\eee
where we applied Lemma \ref{lemBACM}.

\bs
\ni
{\em Proof of Lemma \pref{lemtech2}.} Let $n_0$ be as in Assumption \ref{ass2}, let $\ell \in \NN$ and $u\in C^0([0,T],B^{2n_0+8+\ell })$. Denoting
\be
\label{mup}
u_p=\left\langle u \chi_p\right\rangle,\qquad \mu_p^2=\|u_p \chi_p\|_{C^0([0,T],B^{2n_0+8+\ell})}^2,
\ee
we have
\be
\label{ell15}
\|u\|_{C^0([0,T],B^{2n_0+8+\ell})}^2=\sum_{p\geq 0}\mu_p^2<+\infty.
\ee
From \fref{A2}, we obtain
$$
\|A\left(\frac{t}{\eps^2}\right)u(t)\|_{C^0([0,T],B^{\ell})}\leq C\sum_{p\geq 0}\sum_{q\neq p}(1+q)^{n_0}\,|a_{pq}|\,\|u_q\chi_p\|_{C^0([0,T],B^{\ell+1})}\,,
$$
where we used Assumption \ref{ass2}. Besides, applying Lemma \ref{lemBACM} gives
\bea
&\|u_q\|_{C^0([0,T],H^s(\RR^2))}&=\frac{1}{E_q^{n_0+4+(\ell-s)/2}}\|H_z^{n_0+4+(\ell-s)/2}(I+(-\Delta_{x,y})^{s/2})(u_q\chi_q)\|_{C^0([0,T],L^2)}\nonumber\\
&&\leq C\frac{E_q^{s/2}}{E_q^{n_0+4+\ell/2}}\mu_q
\label{qq}
\eea
for all $s\leq 2n_0+8+\ell$. Hence, from the definition \fref{Bell}, we get
\bee
&\|u_q\chi_p\|_{C^0([0,T],B^{\ell+1})}&\leq CE_p^{(\ell+1)/2}\|u_q\|_{C^0([0,T],L^2(\RR^2))}+C\|u_q\|_{C^0([0,T]H^{\ell+1}(\RR^2))}\\
&&\leq C\frac{E_p^{(\ell+1)/2}+E_q^{(\ell+1)/2}}{E_q^{n_0+4+\ell/2}}\,\mu_q.
\eee
and, by using \fref{borneapq} and \fref{weyl},
$$
(1+q)^{n_0}\,|a_{pq}|\,\|u_q\chi_p\|_{C^0([0,T],B^{\ell+1})}\leq C\frac{E_q^{n_0}}{E_p^2}\,|a_{pq}|\,\frac{E_p^{(\ell+5)/2}+E_q^{(\ell+1)/2}E_p^2}{E_q^{n_0+4+\ell/2}}\,\mu_q \leq C\frac{1}{E_p^2}\,\frac{\mu_q}{E_q}.
$$
Therefore,
$$
\|A\left(\frac{t}{\eps^2}\right)u(t)\|_{C^0([0,T],B^{\ell})}\leq C\left(\sum_{p\geq 0}\frac{1}{E_p^2}\right)\left(\sum_{q\geq 0}\frac{\mu_q}{E_q}\right)\leq C\left(\sum_{p\geq 0}\frac{1}{E_p^2}\right)^{3/2}\left(\sum_{q\geq 0}\mu_q^2\right)^{1/2}
$$
by Cauchy-Schwarz. To conclude, it suffices to use \fref{weyl} and \fref{ell15}: the series converge and we obtain the desired estimate \fref{estiA}.

\bs
\ni
{\em Proof of Lemma \pref{lemtech1}.}  Let $m=4n_0+17$ and let $u\in C^0([0,T],B^m)$ such that $\pa_t u\in C^0([0,T],B^{m-2})$. Denoting now
\be
\label{nup}
u_p=\left\langle u \chi_p\right\rangle,\qquad \nu_p^2=\|u_p \chi_p\|_{C^0([0,T],B^m)}^2+\|\pa_t u_p \chi_p\|_{C^0([0,T],B^{m-2})}^2,
\ee
we have
\be
\label{ell2}
\|u\|_{C^0([0,T],B^m)}^2+\|\pa_t u\|_{C^0([0,T],B^{m-2})}^2=\sum_{p\geq 0}\nu_p^2<+\infty.
\ee
Applying Lemma \ref{lemBACM} as above yields
\be
\label{mm}
E_p^{(m-s)/2}\|u_p\|_{C^0([0,T],H^s(\RR^2))}+E_p^{(m-2-s)/2}\|\pa_tu_p\|_{C^0([0,T],H^s(\RR^2))}\leq C\nu_p
\ee
for all $s\leq m$. By composing the expressions \fref{A2} and \fref{a2} for $A$ and $a$, we obtain
\bee
&A(\tau)a(\tau)u&=i\sum_{p\geq 0}\sum_{q\neq p}\sum_{n\neq q}\frac{e^{i\tau (E_p-E_q)}-1}{E_p-E_q}e^{i\tau (E_q-E_n)}a_{pq}a_{qn}\,\pa^2_x u_n\,\chi_p\\
&&=i\sum_{p\geq 0}\sum_{q\neq p}\frac{1-e^{i\tau (E_q-E_p)}}{E_p-E_q}(a_{pq})^2\,\pa^2_x u_p\,\chi_p\\
&&\quad +i\sum_{p\geq 0}\sum_{q\neq p}\sum_{\scriptsize \begin{array}{l} n\neq q\\ n\neq p\end{array}}\frac{e^{i\tau (E_p-E_n)}-e^{i\tau (E_q-E_n)}}{E_p-E_q}a_{pq}a_{qn}\,\pa^2_x u_n\,\chi_p
\eee
Now, remark that, by \fref{alphap} and \fref{apq}, we have for all $p\in \NN$ the identity
$$1+\sum_{q\neq p}\frac{(a_{pq})^2}{E_p-E_q}=\alpha_p.$$
Therefore we get, using the definition \fref{A0},
\bee
&\left(A(\tau)a(\tau)+i\pa_x^2\right)u&=-iA_0u\\
&&\quad -i\sum_{p\geq 0}\sum_{q\neq p}e^{i\tau (E_q-E_p)}\frac{(a_{pq})^2}{E_p-E_q}\,\pa^2_x u_p\,\chi_p\\
&&\quad +i\sum_{p\geq 0}\sum_{q\neq p}\sum_{\scriptsize \begin{array}{l} n\neq q\\ n\neq p\end{array}}\left(e^{i\tau (E_p-E_n)}-e^{i\tau (E_q-E_n)}\right)\frac{a_{pq}a_{qn}}{E_p-E_q}\pa^2_x u_n\,\chi_p
\eee
and, integrating,
\bea
&&\int_0^t \left(A\left(\frac{s}{\eps^2}\right)+i\pa_x^2\right)a\left(\frac{s}{\eps^2}\right)u(s)ds+i\int_0^tA_0 u(s)ds\nonumber\\
&&\quad =-i\sum_{p\geq 0}\sum_{q\neq p}\frac{(a_{pq})^2}{E_p-E_q}\,\chi_p\int_0^te^{is(E_q-E_p)/\eps^2}\,\pa^2_x u_p(s)\,ds\label{sum}\\
&&\qquad +i\sum_{p\geq 0}\sum_{q\neq p}\sum_{\scriptsize \begin{array}{l} n\neq q\\ n\neq p\end{array}}\frac{a_{pq}a_{qn}}{E_p-E_q}\chi_p\int_0^t\left(e^{is(E_p-E_n)/\eps^2}-e^{i s(E_q-E_n)/\eps^2}\right)\pa^2_x u_n(s)\,ds\nonumber
\eea
In order to estimate the right-hand side of this identity, we claim that, for all $p\in \NN$, $p\in \NN$ and $\lambda\neq 0$, we have
\be
\label{v}
\left\|\chi_p(z)\int_0^te^{i \lambda s/\eps^2}\,\pa^2_x u_q(s,x,y)\,\,ds\right\|_{C^0([0,T],B^1)}\leq C_{T}\frac{\eps^2}{|\lambda|}\,\frac{E_p^{1/2}+E_q^{1/2}}{E_q^{(m-4)/2}}\,\nu_q
\ee
where $C_{T}$ only depends on $T$ and $\nu_n$ is defined by \fref{nup}. This claim is proved below. As a consequence, we can estimate \fref{sum} as follows:\sloppypar
\bee
&&\left\|\int_0^t \left(A\left(\frac{s}{\eps^2}\right)a\left(\frac{s}{\eps^2}\right)+i\pa_x^2\right)u(s)ds+i\int_0^tA_0 u(s)ds\right\|_{C^0([0,T],B^1)}\\
&&\quad \leq C\eps^2\sum_{p\geq 0}\sum_{q\neq p}\frac{(a_{pq})^2}{|E_p-E_q|^2}\,\frac{1}{E_p^{(m-5)/2}}\nu_p\\
&&\qquad +C\eps^2\sum_{p\geq 0}\sum_{q\neq p}\sum_{\scriptsize \begin{array}{l} n\neq q\\ n\neq p\end{array}}\frac{|a_{pq}||a_{qn}|}{|E_p-E_q|}\left(\frac{1}{|E_p-E_n|}+\frac{1}{|E_q-E_n|}\right)\frac{E_p^{1/2}+E_n^{1/2}}{E_n^{(m-4)/2}}\nu_n\\
&&\quad \leq C\eps^2\sum_{p\geq 0}\sum_{q\geq 0}\frac{E_p^{3}}{E_q^2}\,\frac{(1+p)^{2n_0}}{E_p^{(m-5)/2}}\nu_p\\
&&\qquad +C\eps^2\sum_{p\geq 0}\sum_{q\geq 0}\sum_{n\geq 0}(1+q)^{n_0}(1+n)^{n_0}\frac{E_n^{n_0+11/2}}{E_p^2E_q^{n_0+2}}\frac{1}{E_n^{(m-4)/2}}\nu_n\\
&&\quad \leq C\eps^2\sum_{p\geq 0}\sum_{q\geq 0}\frac{1}{(1+q^2)}\frac{\nu_p}{1+p^3}+C\eps^2\sum_{p\geq 0}\sum_{q\geq 0}\sum_{n\geq 0}\frac{1}{(1+p^2)}\frac{1}{(1+q^2)}\frac{\nu_n}{1+n}
\eee
where we used Assumption \ref{ass2}, \fref{borneapq}, \fref{weyl} and recall that $m=4n_0+17$. Hence we deduce \fref{estimI1} by using Cauchy-Schwarz and \fref{ell2}. It remains to prove the claim.

\ms
\ni
{\em Proof of the claim \fref{v}.} Let
\be
\label{vv}
v(t,x,y,z)=\chi_p(z)\int_0^te^{i \lambda s/\eps^2}\,\pa^2_x u_q(s,x,y)\,\,ds,
\ee
for $p\in \NN$, $q\in \NN$ and $\lambda\neq 0$. An integration by parts in \fref{vv} yields
$$
v(t,x,y,z)=i\frac{\eps^2}{\lambda}\chi_p\left(\int_0^te^{i \lambda s/\eps^2}\,\pa^2_x \pa_tu_q(s,x,y)\,\,ds+e^{i \lambda t/\eps^2}\,\pa^2_x u_q(t,x,y)-\pa^2_x u_q(0,x,y)\right).$$
Hence, by using \fref{mm}, we obtain
$$
\|v\|_{C^0([0,T],B^1)}\leq C_{T}\frac{\eps^2}{|\lambda|}\,\frac{E_p^{1/2}+E_q^{1/2}}{E_q^{(m-4)/2}}\,\nu_q,$$
where $C_{T}$ only depends on $T$. This concludes the proof of \fref{v}.

The proof of Lemma \ref{lemtech1} is complete.
\qed

\end{appendix}

\bs

\begin{acknowledgement} The authors were supported by the Agence Nationale de la Recherche, ANR project QUATRAIN. They wish to thank N. Ben Abdallah and F. Castella for fruitful discussions.
\end{acknowledgement}

\end{document}